\documentclass[journal]{IEEEtran}

\IEEEoverridecommandlockouts                              

\usepackage{setspace}
\usepackage{amsmath}
\usepackage{mathalfa}
\usepackage{mathrsfs}
\usepackage{amssymb}
\usepackage{amsthm}
\newtheorem{thm}{Theorem}
\usepackage{color}
\usepackage{subcaption}
\usepackage{soul}
\allowdisplaybreaks
\usepackage{bm}

\usepackage{array}
\usepackage{comment}
\usepackage[ruled,vlined]{algorithm2e}
\newcolumntype{P}[1]{>{\centering\arraybackslash}p{#1}}

\definecolor{orange}{rgb}{1,0.5,0}

\usepackage{subcaption}
\usepackage{graphicx}
\graphicspath{ {pics/} }
\usepackage{flushend}

\begin{document}

\title{\LARGE \bf
Reinforcement Learning versus PDE Backstepping and PI  Control\\for Congested Freeway Traffic
}

\author{Huan Yu,$^{1,*, \dagger}$ Saehong Park,$^{2,\dagger}$ Alexandre Bayen,$^{2}$ Scott Moura,$^{2}$ Miroslav Krstic$^{1}$ 
\thanks{$^{1}$Huan Yu and Miroslav Krstic are with the Department of Mechanical and Aerospace Engineering,
        University of California, San Diego, 9500 Gilman Dr, La Jolla, CA 92093
        {\tt\small email: huy015@ucsd.edu, krstic@ucsd.edu}.}%
\thanks{$^{2}$Saehong Park, Scott Moura and Alexandre Bayen are with the Department of Civil and Environmental Engineering, University of California, Berkeley, CA 94720
        {\tt\small email: sspark@berkeley, smoura@berkeley.edu, bayen@berkeley.edu}.}%
\thanks{$^{*}$ Corresponding author.}
\thanks{$^{ \dagger}$ Huan Yu and Saehong Park contributed equally to this work.}
}

\maketitle
	
	
	
\IEEEpeerreviewmaketitle

\begin{abstract}
		We develop reinforcement learning (RL) boundary controllers to mitigate stop-and-go traffic congestion on a freeway segment. The traffic dynamics of the freeway segment are governed by a macroscopic Aw-Rascle-Zhang (ARZ) model, consisting of $2\times 2$ quasi-linear partial differential equations (PDEs) for traffic density and velocity. Boundary stabilization of the linearized ARZ PDE model has been solved by PDE backstepping, guaranteeing spatial $L^2$ norm regulation of the traffic state to uniform density and velocity and ensuring that traffic oscillations are suppressed. Collocated Proportional (P) and Proportional-Integral (PI) controllers also provide stability guarantees under certain restricted conditions, and are always applicable as model-free control options through gain tuning by trail and error, or by model-free optimization. Although these approaches are mathematically elegant, the stabilization result only holds locally and is usually affected by the change of model parameters. Therefore, we reformulate the PDE boundary control problem as a RL problem that pursues stabilization without knowing the system dynamics, simply by observing the state values. The proximal policy optimization, a neural network-based policy gradient algorithm, is employed to obtain RL controllers by interacting with a numerical simulator of the ARZ PDE. Being stabilization-inspired, the RL state-feedback boundary controllers are compared and evaluated against the rigorously stabilizing controllers in two cases: (i) in a system with perfect knowledge of the traffic flow dynamics, and then (ii) in one with only partial knowledge. We obtain RL controllers that nearly recover the performance of the backstepping, P, and PI controllers with perfect knowledge and outperform them in some cases with partial knowledge. It must be noted, however, that the RL controllers are obtained by conducting about one thousand episodes of iterative training on a simulation model. This training cannot be performed in a collision-free fashion in real traffic, nor convergence guaranteed when training. Thus, we demonstrate that RL approach has learning (i.e. adaptation) potential for traffic PDE system under uncertain and changing conditions, but RL is neither simple nor a fully safe substitute for model-based control in real traffic systems. 

\end{abstract}

\section{Introduction}

 Traffic congestion has been an inescapable problem in large and growing metropolitan areas across the world. Among various traffic congestion phenomena, stop-and-go traffic~\cite{bell}~\cite{Daganzo:11} ~\cite{Seibold1}~\cite{Seibold2} is a common variety that appears every day in congested freeways, causing unsafe driving conditions, more fuel consumption and emissions. Traffic congestion is characterized by the propagation of oscillatory waves on roads and caused by delayed driver response.
 Freeway traffic management usually relies on static road infrastructure to regulate traffic flow, such as ramp metering and varying speed limits (VSL).

Depending on the employed traffic model and control objectives, numerous model-based traffic control strategies have been developed using Lyapunov stability analysis. Macroscopic traffic modeling is particularly well-suited to control design since it describes the overall spatial-temporal dynamics of traffic flow on a freeway segment. The state variables of the model are aggregated traffic values. They are continuous in time and space, easy to be measured and to be regulated in transportation management. The first-order hyperbolic partial differential equation model by Lighthill, Whitham and Richards (LWR)~\cite{LW}~\cite{R} is a conservation law of traffic density, describing the formation and propagation of traffic density waves on the road. The oscillatory, unstable behaviors observed in stop-and-go traffic require a higher-order model rather than the static LWR model. The second-order Aw-Rascle-Zhang (ARZ) model~\cite{AW}~\cite{Zhang} addressed this issue by adding another PDE for the velocity state, leading to a $2\times 2$ nonlinear coupled hyperbolic PDE system. The authors in \cite{Seibold_data} compared the data-fitted LWR and ARZ model and found out that the second-order ARZ model is significantly more accurate in predicting traffic velocity. Therefore, we examine the state-of-art ARZ model to capture stop-and-go traffic oscillations in this paper, and study the boundary control problem.

\subsection{Lyapunov-based control approach and challenges}
The control of traffic PDE models using VSL and ramp metering has been addressed by two categories of research. The first category reduces the infinite-dimensional traffic PDE model to a system of finite-dimensional ordinary differential equations, and then applies optimal control approaches. The authors in \cite{Gomes:06} proposed a cell-transmission model (CTM) derived from the LWR model. The authors of \cite{Papageorgiou:90} proposed the ``METANET'' model which is discretized from the second-order traffic model. Many recent efforts ~\cite{Bekiaris-Liberis:18}~\cite{bell}~\cite{Karafyllis:18}~\cite{Yu:auto19}~\cite{Huan2}~\cite{Yu:TCST}~\cite{Yu:TAC}~\cite{Yu:JDSMC}~\cite{Liguo:19} focused on developing feedback control strategies directly for the macroscopic traffic PDE models, instead of their approximated versions. The Lyapunov approach was applied for stability analysis. In a nutshell, the Lyapunov-based control approach for traffic PDE models consists of two steps: (i) obtain the state-feedback controllers through various control design methods, and (ii) use Lyapunov analysis to study the closed-loop system dynamics.
 
Studies that examine the first-order LWR model include ~\cite{Bekiaris-Liberis:18}~\cite{Karafyllis:19}~\cite{Yu:TAC}~\cite{Yu:JDSMC}, and all investigate the traffic bottleneck problem. Predictor feedback control laws were designed to compensate the delay PDE, i.e, the linearized LWR model ~\cite{Bekiaris-Liberis:18}~\cite{Yu:TAC}. Extremum seeking, a non-model based real-time adaptive control method was used in~\cite{Yu:JDSMC} to find the unknown optimal density of the bottleneck area. In \cite{Karafyllis:19}, stabilization of a spatially uniform equilibrium profile for the LWR model was obtained under in-domain VSL control. All of these studies achieved closed-loop stability using Lyapunov analysis.

Studies that examine the second-order ARZ model include~\cite{Karafyllis:18}~\cite{Yu:auto19}~\cite{Liguo:19}, and this paper. In \cite{Karafyllis:18}, a nonlinear boundary feedback law was designed that controls the inlet flow and achieves global stabilization for a modified second-order ARZ model. PI controllers were designed in~\cite{Liguo:19} for inlet traffic flow rate control through ramp metering and outlet velocity control through VSL. In~\cite{Yu:auto19}, a backstepping controller was developed for ramp metering at the outlet, to regulate upstream traffic. A collocated P controller was also proposed for inlet flow control of the downstream traffic. The aforementioned control designs are representatives of boundary control of hyperbolic PDE systems, and guarantee Lyapunov stability of the closed-loop system. 

Despite differences between the control design approaches, the key idea behind each method is similar. Namely, each method analyzes the particular structure of the PDE system to design a boundary feedback control law and then uses Lyapunov analysis to obtain closed-loop stability. In~\cite{Liguo:19}, the authors analyzed PI boundary controllers for the linearized and homogeneous ARZ model and exponential stability of the closed-loop system was obtained using Lyapunov stability analysis. A boundary controller for the linearized homogeneous ARZ PDE is developed in~\cite{Yu:auto19}~\cite{Huan2}. A linear backstepping state-feedback controller and a static boundary feedback controller were proposed to achieve $L^2$ norm stabilization of traffic oscillations in finite time.

Despite the theoretical results from the Lyapunov-based controllers, several challenges still exist and include:

\begin{itemize}
    \item {\bf Nonlinearity:} The backstepping, P and PI boundary control methods discussed above are developed for linearized PDEs. Thus only local stabilization is guaranteed for the nonlinear system. The linearization approach may fail for large perturbations from the linearization point.

    \item {\bf Model knowledge}: Both model-based designs (e.g. PDE backstepping) and model-free designs (e.g. PI control) involve gain selection subject to certain theoretically-determined bounds, which require some knowledge of the system such as real-time full state measurements as well as model parameters that characterize the overall system dynamics. The parameter values are usually obtained by calibrating the PDE model with traffic field data. This model calibration process can be laborious, and specific to certain traffic conditions. In addition, the deterministic traffic PDE model cannot address random events such as traffic accidents, changes in road attributes, etc. 
    
    \item {\bf Adaptive optimization:} Besides traffic flow stabilization, there are other relevant performance metrics, including fuel consumption, travel time, etc. When another objective is considered, then one needs to re-formulate and re-design the controller with a Lyapunov-based approach. Different design tools might be required. It is desirable to have an approach with modest tuning that can adapt to various problems and objectives arising in traffic systems.

\end{itemize}
The challenges mentioned above can be addressed by reinforcement learning (RL). This model-free approach relies on few assumptions about the system, including being agnostic about nonlinearities. From a control law perspective, RL produces a neural network-based feedback controller that maps observations (states) into actions (control inputs) through training the networks. A RL formulation to the traffic control problem has the potential to enable an \textit{expert-free} traffic management system, in the sense that it does not require expert knowledge of the system dynamics or a control-theoretic training. Neural networks are used as parametric approximators for functions of states whose weights are updated automatically to learn the optimal control action. The model complexity and identification process faced by Lyapunov-based controllers are circumvented. RL controllers are trained to maximize the reward function, which in this paper is chosen as the $L^2$ spatial norm of states to achieve stabilization. However, the proposed formulation can be easily adapted for different objectives by changing the reward formulation.

Extremum Seeking, a real-time adaptive optimization technique, is an alternative to RL for tuning controller gains in a model-free fashion. The algorithms in~\cite{Frihauf}~\cite{Nesic}~\cite{killingsworth} can be applied for offline learning of feedback laws based on simulation models, as with RL and the ARZ model in this paper. In contrast, the algorithm in~\cite{Radenkovic} is designed for online tuning of feedback laws. In extremum seeking, a small excitation is used to perturb the feedback gains being tuned and to produce estimates of the gradient of a cost function. Convergence to a neighborhood of the optimal gain values is proved by means of averaging analysis and singular perturbation theory. Similarly, the adaptive PDE backstepping control design in \cite{Yu:auto19} achieves output feedback stabilization with gradient-based estimators for some unknown model parameters. In this paper, we employ the RL boundary control approach to deal with the nonlinear ARZ model with uncertain steady state conditions, which has not been studied before.

\subsection{RL boundary control approach}
Recent developments in Reinforcement Learning (RL) have enabled model-free control of high-dimensional continuous control systems. Model-free RL control does not require any prior knowledge of the model structure, nor model calibration. Instead, the control input is obtained by iterative interactions with {\it a simulator of the system dynamics}. Thus the PDE boundary control problem is formulated as a RL problem, requiring minimal assumptions on the infinite dimensional system yet in practice obtains reasonably good controllers through a completely data-driven process. 

Learning-based methods have drawn attention from both the PDE modeling community and transportation researchers. Algorithms for solving high-dimensional parabolic PDE systems were derived using deep learning in \cite{JHan18} by reformulating the PDEs as stochastic differential equations and approximating the gradient of unknown solutions with deep neural networks (DNN). More recent work in \cite{YWang19} proposed RL-learned numerical solutions for scalar conservation laws which autonomously generate accurate numerical schemes for various situations. For the PDE control problem, boundary control of time-varying 2D convection-diffusion PDEs with an application for heating, ventilating, air conditioning (HVAC) control design was developed in \cite{Farahmand17,Farahmand18}. 

In the domain of traffic management, researchers have been applying RL to various traffic problems. Different levels of traffic modeling are employed, depending on the problem. For example, the authors in \cite{LiLi} examined traffic light signal timing at one intersection using deep Q-learning with a traffic simulator, Paramics, which is based on microscopic traffic models. Studies on the traffic signal scheduling problem were further extended to multi-agent control for five centrally-connected traffic intersections~\cite{IArel}, which was described with a traffic queuing model. The article~ \cite{MARLIN} used a multi-agent RL algorithm to control the traffic light around a traffic junction. The authors proposed a framework where each agent was able to switch between independent and integrated modes in which the agent solved the multi-agent RL problem using modular Q-learning. The authors of~\cite{Wu17} formulated a deep RL framework for mixed-autonomy traffic in some experimental scenarios using the traffic microsimulator SUMO. Under the same framework, \cite{Vinitsky18} developed RL controllers for connected autonomous vehicles to de-congest traffic bottlenecks.

There are few RL works related to macroscopic traffic models in the literature, although the PDE model is particularly well-suited for modeling congested traffic flow patterns. The reason is that RL control of PDEs involves high-dimensional state spaces which makes the approximation of value functions challenging. The authors in \cite{BELLETTI} considered the cell transmission model which is obtained from discretization of the LWR first-order PDE model. A RL-based controller was designed using different policy gradient methods, such as REINFORCE, Trust Policy Optimization (TRPO), and the Truncated Natural Policy Gradient (TNPG) algorithm. Incoming traffic flow is actuated by the RL controller such that the traffic flow is optimized for some target outflow. Our work on RL control differs on the traffic problem to be solved, the model being used, and the methodology being employed, which guides our study and analysis from a very different perspective. We focus on the stabilization problem, and examine a second-order PDE model to describe stop-and-go traffic oscillations. The resulting RL controllers are tested, and we examine under what circumstances RL could be a better choice relative to a Lyapunov-based PDE controller.

\subsection{Contribution}
We formulate a state regulation control problem for the ARZ PDE model via boundary control. Then we develop a RL control approach based on the PPO algorithm, which falls within the class of policy gradient methods. PPO ultimately yields a state-feedback boundary controller from iterative interactions with a simulation environment, as opposed to direct synthesis from a mathematical model. The performance of the RL controller is compared with a PDE backstepping controller, a P controller and PI controllers. 

Interestingly, RL controllers nearly recover the stabilization performance of the Lyapunov-based PDE control approaches for a system with perfect knowledge of the model. However, in a system with partial knowledge where the steady state traffic is lighter or denser than what we assumed in constructing Lyapunov-based controllers, the RL controller obtained from a stochastic training process outperforms the Lyapunov-based controllers.
However, it must be noted that the RL controllers are obtained by conducting about one thousand training episodes on a simulation model. Collision-free training is not possible in real life traffic, nor is the iterative training guaranteed to converge. Although RL demonstrates learning (i.e. adaptation) potential under uncertain and changing conditions, it is neither simple nor a fully safe substitute for model-based control in real traffic systems.

The contributions of this article are twofold: First, although several Lyapunov-based boundary controllers have been designed for the linearized ARZ model, there is no result for control of the original nonlinear ARZ model. This motivates the RL control in this paper, which provides the first native controller for a nonlinear ARZ PDE model, to the authors' best knowledge. Second, the performance of the RL controllers are compared with PDE backstepping, P, and PI controllers for stabilization of stop-and-go traffic, and are assessed for different scenarios. By evaluating performance with full knowledge of the dynamics and one with partial knowledge, we demonstrate via several numerical simulations under what circumstances RL could be a better choice from a traffic application perspective.   


The outline of this article is as follows: Section 2 summarizes the ARZ traffic flow model and the Lyapunov-based boundary controllers. Section 3 details the reinforcement learning approach for boundary control of the ARZ PDE model. Section 4 includes the comparative analysis and numerical results, comparing the closed-loop performance of setpoint boundary inputs, backstepping, P, and PI controllers in different scenarios. The conclusion summarizes the main results, and discusses future work.

\section{Problem Statement}
\subsection{ARZ PDE Traffic Model}
We consider the ARZ PDE model to describe the traffic dynamics on a freeway segment~\cite{AW}~\cite{Zhang}. The state variables are traffic density $\rho(x,t)$ and traffic speed $v(x,t)$, defined on the domains $x\in[0,L]$, $t\in[0,T]$:
\begin{align}
\rho_t + (\rho v)_x =& 0,\label{rho} \\
(v-V(\rho))_t + v(v-V(\rho))_x =& \frac{V(\rho)-v}{\tau}, \label{v}
\end{align}
where $(\cdot)_z$ is short-hand notation for the differential operator $\partial/\partial z$. Parameter $\tau$ is the relaxation time, and captures how quickly drivers adjust their velocity to the equilibrium.
The equilibrium velocity-density relationship $V(\rho)$ is a decreasing function, e.g., the linear Greenshield's model 
$V(\rho) = v_m \left(1 - \frac{\rho}{\rho_m}\right)$
representing the reduction of velocity with the increase of density. The traffic flow rate is defined with 
\begin{align}
    q = \rho v.
\end{align}
Our control objective is to regulate the state around an equilibrium reference state $(\rho^\star, v^\star)$, where
\begin{align}\label{ref-state}
v^\star = V(\rho^\star),
\end{align}
satisfies the equilibrium density and velocity relation. We choose the density $\rho^\star$ such that the reference system $(\rho^\star, v^\star)$ is in the congested regime for dense traffic, which can be characterized by the two characteristics of the linearized PDE model \cite{Yu:auto19}. 
\begin{align}
\lambda_1 = &v^\star>0, \label{char1} \\ 
\lambda_2 =& v^\star + \rho^\star V'(\rho^\star)<0. \label{char2}
\end{align}
The first characteristic is always greater than $0$. When traffic is light, $\lambda_2 >0$ is satisfied. When traffic is dense, then $\lambda_2 < 0$ and in this regime there can be upstream propagation of oscillations in the states. This can also be characterized with the Traffic Froude Number (TFN) in \cite{BELLETTI}. Consequently, we have hetero-directional propagation of oscillations in congested traffic. The density oscillations are carried downstream by vehicles while the velocity oscillations are transported upstream. Intuitively, drivers are mostly affected by vehicles driving in front of them. The stop-and-go traffic is characterized by oscillations, caused by delayed driver reaction to vehicles in front of them. 
\begin{figure}
	\centering
\includegraphics[width=6.5cm]{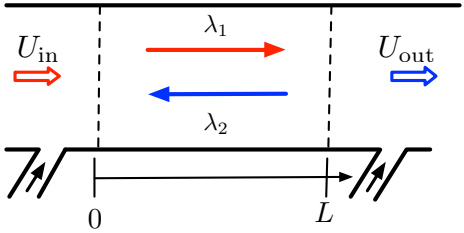}
	\caption{Stabilization of stop-and-go traffic on freeway with ramp metering located at boundaries of the segment. }
	\label{ARZ_schematic}
\end{figure}
As shown in Fig.~\ref{ARZ_schematic}, the control inputs $U_{\rm in}$ and $U_{\rm out}$ are transported from the actuated boundaries to in-domain states with the two characteristic speeds respectively.

\subsection{Boundary control design}
Our control objective is $L^2$ regularization of the density and velocity states to uniform steady state values by ramp-metering on the boundary. The oscillations can occur due to delayed driver response without ramp-metering control.
In order to reduce oscillations in congested traffic, we actuate traffic flow from the upstream inlet and downstream outlet of the freeway segment, using on-ramp metering as shown in Fig.~\ref{ARZ_schematic}. Alternatively, one can actuate velocity by installing VSLs.

Boundary control algorithms are designed to stabilize the traffic around the reference steady state. For the ARZ PDE model, the control design guarantees that the state variables $(\rho(x,t),v(x,t))$ are regulated to the reference system $(\rho^\star, v^\star)$ in the spatial $L^2$ norm, i.e. 
\begin{align}
&||\rho(x,t) - \rho^\star|| \to 0,\\
&||v(x,t) - v^\star|| \to 0,
\end{align}
where $||\cdot||$ is the spatial $L^2$ norm, i.e., $||f(x,t)|| = \left(\int_{0}^{  L}f^2(x,t)dx\right)^{1/2}$.

For the segment of freeway traffic, boundary actuation is implemented with ramp metering controlling the traffic flow rate entering from the on-ramp to the mainline. Both inlet and outlet flow rates can be controlled through ramp metering
\begin{align}
  \rho(0,t) v(0,t) =&  U_{\rm in} (t),\\
  \rho(L,t) v(L,t) =&  U_{\rm out} (t),
\end{align}
We denote the boundary control inputs of different strategies with $U_{\rm in,\star} (t), U_{\rm out, \star} (t),$ where $\star$ represents the name of the control design. As summarized in Table~\ref{bcd}, the Lyapunov-based boundary control algorithms consist of setpoint boundary inputs, backstepping PDE control~\cite{Yu:auto19}, Proportional~\cite{Yu:auto19} and Proportional Integral control~\cite{Liguo:19}. Despite different actuation locations and assumptions, three different state-feedback controllers achieve Lyapunov stabilization of the PDE system by collectively shifting the linearized PDE eigenvalues onto the open left half plane. The control gains used for the backstepping method are computed numerically~\cite{Yu:auto19}. The static proportional gain in inlet control is derived by solving the linearized PDE system analytically~\cite{Yu:auto19}. The static proportional and integral gains in PI control are obtained through trial and error.

\begin{table}[t!]
	\caption{Lyapunov-based Boundary Control Design}
	\label{bcd}
	\centering
	\begin{tabular}{P{3cm}|P{2cm}|P{2cm}}
		\hline \hline
		\textbf{Design method} & $U_{\rm in}$ & $U_{\rm out}$ \\
		\hline
		 Setpoint & $q^\star$ & $q^\star$  \\
     	 Backstepping  & $q^\star$ & $U_{\rm out, BKST}$   \\
		Proportional  & $U_{\rm in, L}$ & $q^\star$    \\
	 	Proportional Integral  & $U_{\rm in, PI}$ & $U_{\rm out, PI}$   \\
		\hline \hline
	\end{tabular}
\end{table}
\subsubsection{Setpoint control}
We simply limit the outgoing flow rate to be the setpoint state value as the boundary control. This simple method is considered as a baseline case to compare against more sophisticated control designs. The constant incoming and outgoing flow rate at the boundaries are implemented as,
\begin{align}
U_{\rm in, O} (t) =& q^\star,\label{bc1}\\
U_{\rm out, O} (t) =& q^\star.\label{bc2}
\end{align}
where $q^\star = \rho^\star v^\star $ is the reference steady state flow rate. The ARZ PDE system with this setpoint control can be unstable (in the linear sense) for certain combinations of $q^\star$ and model parameters $V(\rho)$ and $\tau$~\cite[Section III]{Huan2}.

\subsubsection{PDE backstepping control} \label{bkst_method}
The backstepping controller produces a full-state-feedback control law located at the outlet boundary. The control design is developed in \cite{Yu:auto19} for the linearized ARZ PDE model, and is applied to the original nonlinear ARZ model. For initial conditions near the reference, the system is locally exponentially stable under outlet actuation. The backstepping method is used to construct the full-state-feedback control law \cite{intro_bkst}. The key idea of the backstepping method is to transform the original linearized ARZ model to a target system where the undesirable instabilities in the PDEs are transformed to the outlet boundary. The controller is designed such that the instability is cancelled out and an exponentially stable target system with zero boundary input is obtained for the closed-loop system. 

The boundary condition at the inlet, and the full-state-feedback outlet backstepping controller are given by
\begin{align}
U_{\rm in, BKST} (t) =& q^\star,\label{bkst_in}	\\
\notag U_{\rm out, BKST}(t) =&q^\star + \rho^\star \int_{0}^{L} c_v(\xi)(v(\xi,t) -v^\star)d\xi \label{bkst_out} \\ 
 & + \int_{0}^{L}c_q(\xi) (q(\xi,t)-q^\star) d\xi,
\end{align}
where the control gains are given by
$
c_v(\xi) = M(L-\xi) + \frac{\lambda_2}{\lambda_1}K(L,\xi)\exp\left(\frac{ \xi }{\tau v^\star }\right),
c_q(\xi) =\frac{\lambda_1 -\lambda_2}{\lambda_1}K(L,\xi)\exp\left(\frac{ \xi }{\tau v^\star} \right)
$
for $\xi \in [0,L]$. The control gain kernels $K(L,\xi)$, $M(L-\xi)$ are obtained by solving hyperbolic equations that govern the kernel variables $K(x,\xi)$ in a triangular spatial domain $\mathcal{T}=\{(x,\xi): 0 \leq \xi \leq x \leq 1\}$. The numerical solution of the kernel equations is easily obtained, given the model parameters and steady states. The equilibrium $(\rho^\star,  v^\star)$ of the linearized system is exponentially stable in the $L^2$ sense and the equilibrium is reached in a finite time.
\begin{thm}[{\cite[Theorem 2]{Yu:auto19}}]
	Consider system \eqref{rho}-\eqref{v} 
	linearized around \eqref{ref-state}, initial conditions  $(\rho(x,0),  v(x,0)) \in L^2([0,L])$ 
	and the control laws \eqref{bkst_in}, \eqref{bkst_out}. The equilibrium $(\rho^\star,  v^\star)$ of the linearized system is exponentially stable in the $L^2$ sense and the equilibrium is reached in a finite time $t= \frac{L}{|\lambda_1|} + \frac{L}{|\lambda_2|}$.
\end{thm}
The proof of Theorem 1 is detailed in \cite{Yu:auto19}, using Lyapunov analysis.

\subsubsection{P control} \label{p_method}
The proportional controller is an output-feedback control law that actuates traffic flow rate at the inlet boundary, as proposed in \cite{Yu:auto19}. The control input only requires online measurements of velocity at the inlet boundary, collocated with the actuation location.  The key intuition of the control design is to cancel the forward coupling in the system and then the closed-loop system can be directly solved. Local exponential stability and finite-time convergence is guaranteed. The boundary condition at the outlet is a constant flow rate and the inlet boundary is actuated by the collocated output-feedback controller,
\begin{align}
U_{\rm in, L}(t) =&q^\star + g_{P}({v}(0,t)-v^\star),\label{p_in}\\
U_{\rm out, L}(t) =& {q}^\star,\label{p_out}
\end{align}
where $g_P = \rho^\star + {v^\star}/{ V'(\rho^\star)}$ is a constant control gain.

\begin{thm}[{\cite[Theorem 1]{Yu:auto19}}]
	Consider system \eqref{rho}-\eqref{v} 
	with initial conditions  $(\rho(x,0),  v(x,0)) \in L^2([0,L])$,
 linearized around \eqref{ref-state} 
	and the control law \eqref{p_in}, \eqref{p_out}. The equilibrium $(\rho^\star,  v^\star)$ of the linearized system is exponentially stable in the $L^2$ sense and the equilibrium is reached in a finite time $t= \frac{L}{|\lambda_1|} + \frac{L}{|\lambda_2|}$.
\end{thm}

\subsubsection{PI control} \label{PI_method}
In the proportional-integral controller, feedback control is applied to both the inlet and outlet boundary values. The PI boundary feedback controllers are proposed in \cite{Liguo:19} for the linearized ARZ model, which guarantees local exponential stability of the closed-loop system. The exit velocity at the outlet is controlled such that $v(L,t) = U_{\rm out, PI}(t)$. Ramp metering controls the incoming flow rates, given by
\begin{align}
U_{\rm in, PI}(t) =&q^\star + k_{P}^r({\rho}(L,t)-\rho^\star) + k_I^r  \int_{0}^{t}({\rho}(L,t)-\rho^\star)  ds, \label{PI_in}\\
U_{\rm out, PI}(t) =&{v}^\star + k_{P}^v({v}(0,t) - v^\star) + k_I^v  \int_{0}^{t}({v}(0,t) - v^\star)ds,\label{PI_out}
\end{align}
where $k_{P}^v , k_I^v$, $k_{P}^r$, $k_I^r$ are tuning gains. For this anti-collocated output-feedback structure, a set of linear matrix inequalities are given for allowable control gains that guarantee Lyapunov stability. Within these conditions, the specific values are obtained through trial and error, as detailed in \cite{Liguo:19}.

\begin{thm}[{\cite[Theorem 1]{Liguo:19}}]
	Consider system \eqref{rho}-\eqref{v} 
	linearized around \eqref{ref-state} 
	and the control laws \eqref{PI_in}, \eqref{PI_out}. The equilibrium $\rho(x,t) \equiv \rho^\star,  v(x,t) \equiv v^\star$ of the linearized system is exponentially stable in the $L^2$ sense.
\end{thm} 

\begin{algorithm}[t!]
\SetAlgoLined
\KwResult{Feedback boundary stabilization}
{\bf initial conditions:} $\bm s_0 = \rho(x,0), v(x,0)$; 

{\bf Input}:  parameters $V(\rho), \tau$, steady state $\bm s^\star = (\rho^\star, v^\star)$\;

 \For{$t=1$ \KwTo $N-1$}{
 Compute Lyapunov-based boundary controllers $U_{\rm in,out} = u \left(\bm s_t, V(\rho), \tau, \bm s^\star\right)$, where the functional of boundary control $u$ depends on the choice of control design in Table.1
 \\

Update the state variables $\bm s_{t+1}$ with actuated boundary control inputs $U_{\rm in,out}$ and $\bm s_{t}$\;}

\caption{Boundary Control Scheme}
\label{Boundary Control Algorithms}
\end{algorithm}

All of the aforementioned controllers are designed from the linearized ARZ model, and therefore the stabilization guarantees are local. In addition, model uncertainty has not been taken into consideration, such as random events that frequently appear in traffic. This motivates us to explore a RL approach.

\section{Control of ARZ model with Reinforcement Learning}
In this section, we introduce a reinforcement learning approach for boundary control of the nonlinear ARZ traffic flow model. Although explicit knowledge of the PDE model is not required, RL considers that the traffic dynamics are governed by a Markov Decision Process (MDP). In particular, we will use a policy gradient method since they are applicable to continuously valued control actions.

\begin{algorithm}[t!]
\SetAlgoLined
 Initialize parameters for actor $\theta_{0}$ and critic $\phi_{0}$ \\
 \For{$ k = 1$ \KwTo {$E$} }{
  Initialize the states $s_{0}$ from random distribution \\
  \For{$t=1$ \KwTo $N$}{
   Compute reward $r_{t}$  = $\lVert s_{t} - {s^\star} \rVert_{2}^{2}$; \\
  Draw an action $a_{t}$ from stochastic policy, i.e., $a_{t} \sim \mathcal{N}(\mu(s_{t}),\sigma^{2}(s_{t}))$; \\
  $a_{n}$ goes to boundary conditions of PDE system; \\
  Update state $s_{t+1}$ using $a_t$ and $s_t$ ; \\
  }
  Collect set of trajectories $\mathcal{D}$ driven by policy $\pi_{\theta_{k}}$ \\
  Compute total discounted reward $R_{t}$. \\
  Compute advantage estimates, $A_{t}$ from critic $V_{\phi_{k}}$. \\
  Update the actor, $\theta_{k+1}$ by \eqref{eqn:actor_update}. \\
  Update the critic, $\phi_{k+1}$ by regression.
 }
 \caption{RL Control Procedure}
 \label{algorithm:RL_Control}
\end{algorithm}

\subsection{Boundary control of PDE as a MDP}

In the MDP setting, we seek the policy that maximizes the total reward received from the environment, i.e. the plant. At each time step $ t$ the environment conditions are described by a state vector, $\bm{s}_t\in \mathcal{S}$, where $\mathcal{S}$ is the state space while the control policy picks an action $\bm{a}_{t}\in \mathcal{A}$, with $\mathcal{A}$ being the action space. The control policy is a full state feedback control law that selects an action $\bm{a}_{t}$ based on an observation of the state $\bm{s}_{t}$. The action is applied to the environment, whose state evolves to $\bm{s}_{t+1}\in \mathcal{S}$, according to the state-transition  probability  $\mathcal{P}(\bm{s}_{t+1}|\bm{s}_{t},\,\bm{a}_{t})$, and the agent receives a scalar reward $\bm{r}_{t+1}=r(\bm{s}_t,\bm{a}_t)$. The policy is represented by $\pi$ which maps the state to the action and can be either deterministic or stochastic. The total discounted reward from time $t$ onward can be expressed as:
\begin{align}
	R_{t} & =\sum_{k=0}^{\infty} \gamma^{k}r(\bm{s}_{t+k},\bm{a}_{t+k}) \label{eqn:reward_fct}
\end{align}
where $\gamma \in [0,1]$ is the discount factor. Note, we have abused notation by allowing $t$ to represent both continuous and discrete-valued times. Nevertheless, the meaning will be clear from context. 

The implementation of the RL controller is summarized in Algorithm~\ref{Boundary Control Algorithms}. We represent the state at time $t$ with $\bm s_t = (\rho(\cdot,t), v(\cdot,t))$. When the nonlinear ARZ model is boundary actuated, then it naturally forms a sequential decision making problem, which is modeled as a MDP. 
 
We consider a discretized approximation of the nonlinear ARZ PDE model using the second-order Lax-Wendroff scheme \cite{LaxWendroff} with conservative state variables. The solution $\rho(x,t)$ and $v(x,t)$ to the ARZ PDE model is approximated by piecewise constant functions on discretized temporal and spatial domains. The solution domain is $[0, L] \times [0, T]$. The discretization resolution $\Delta t =  T/(N-1)$ and $\Delta x= L/(M-1)$ are chosen such that the Courant-Friedrichs-Lewy (CFL) condition is met, i.e. $\Delta t \leq c\Delta x$, where $M, N$ are the number of nodes for the spatial and temporal domains respectively and $c$ is defined as the maximum characteristic speed of the nonlinear hyperbolic ARZ PDE model $\frac{\Delta x}{\Delta t} \geq \max |\lambda_{1,2}|$. We can compactly write the discretized nonlinear ARZ PDE system  \eqref{rho} -- \eqref{v} as
 \begin{align}
 \bm s_{t+1} = f(\bm s_t, u_t),
 \end{align}
where $\bm s_t$ describes the density and velocity state at t, $u_t$ represents the boundary control inputs, depending on the choice of control design in Table 1. The function $f$ represents the discretized deterministic dynamics for the temporal evolution of the PDE system. The PDE dynamics at the current time instant are fully described given the current state $s_t$ and control inputs $u_t$. The deterministic temporal evolution of the PDE system can be generalized to stochastic dynamics by the MDP state transition probability
 \begin{align}
     \bm s_{t+1} \mathtt{\sim} \mathcal{P}(\bm{s}_{t+1}|\bm{s}_{t},\,\bm{a}_{t}).
 \end{align}

The discretized states, $\bm{s}_{t}$, and boundary control input, $\bm{a}_{t}$ are written as:
\begin{align}
 \notag \bm{s}_{t} =& [\rho(0,t), \rho(\Delta x,t), \cdots, \rho(L,t),\\ &\; v(0,t), v(\Delta x, t), \cdots, v(L,t)]^{\top},\\
 \bm{a}_{t} =& [q(0,t),q(L,t)]^{\top}.
\end{align}
where $\rho(\cdot,t)$ and $v(\cdot,t)$ are the traffic density and velocity that are discretized in the spatial and temporal domains.

The reward $r(\bm{s},\bm{a})$ is defined by the $L^2$ norm of the states, namely,
\begin{align}\label{reward}
 r_t(\bm{s}_{t},\bm{a}_{t}) =& - \left[\frac{ \Sigma_i \ \rho(i\cdot \Delta x,t) - \rho^\star}{\rho^\star}\right]^{2} \nonumber \\
 & - \left[ \frac{ \Sigma_i v(i \cdot \Delta x,t) - v^\star}{v^\star} \right]^{2}.
\end{align}
The reward is equivalent to the control objective that achieves regulation of the traffic states to a spatially uniform density and velocity.

\subsection{Reinforcement Learning Background}
In the following, we briefly review some essential reinforcement learning concepts to aid the subsequent discussion. A thorough exposition can be found in \cite{bertsekas2005dynamic}~\cite{sutton2000policy}. 

The \textit{state value function}, $V^{\pi}(\bm{s}_t)$ is the expected total discounted reward starting from state $\bm{s}_t$. In the controls community, this is sometimes called the cost-to-go, or reward-to-go. Importantly, note that the value function depends on the control policy. If the agent uses a given policy $\pi$ to select actions starting from the state $\bm{s}_t$, then the corresponding value function is given by:
\begin{align}
	V^{\pi}(\bm{s}_t) = \mathbb{E}\Big[ R_{t} \ | \  \bm{s}_{t} \Big]
    \label{eqn:Vfct}
\end{align}

Then, the optimal policy $\pi^{*}$ is the policy that corresponds to the maximum  value $V^{{*}}(\bm{s}_t)$ of the value function
\begin{align}\label{eqn:Vmax}
	\pi^{*} = \text{arg}\max_{\pi} V^{\pi}(\bm{s}_t)
\end{align}
The solution of \eqref{eqn:Vmax} is pursued by Dynamic Programming (DP) methods. DP, however, requires knowledge of the model / environment.

The next definition, known as the ``Q-function'', is fundamental since it enables the concept of model-free reinforcement learning.  Consider the \textit{state-action value function}, $Q^{\pi}(\bm{s}_t,\bm{a}_t)$, which is a function of the state-action pair that returns a real value. In other words, it corresponds to the expected total discounted reward when the action $\bm{a}_t$ is taken in state $\bm{s}_t$, and then the policy $\pi$ is followed henceforth. Mathematically,
\begin{align}
    Q^{\pi}(\bm{s}_t,\bm{a}_t)  = \mathbb{E}\Big[ R_{t} \ | \ \bm{s}_{t},\, \bm{a}_{t} \Big] \label{eqn:Qvalue}
\end{align}

The optimal Q-function is given by
\begin{align}
Q^{*}(\bm{s}_t,\bm{a}_t)=\text{arg}\max_{\pi} Q^{\pi}(\bm{s}_t,\bm{a}_t)
\end{align}
and represents the expected total discounted reward received by an agent that starts in $\bm{s}_t$, picks (possibly non-optimal) action $\bm{a}_t$, and then behaves optimally afterwards. Since $V^{{*}}(\bm{s}_t)$ is the maximum expected total discounted reward starting from state $\bm{s}_t$, it will also be the maximum of $Q^{*}(\bm{s}_t,\bm{a}_t)$ over all possible actions $\bm{a}_t\in \mathcal{A}$,
\begin{align}
	V^{*}(\bm{s}_t) = \max_{\bm{a}_t\in \mathcal{A}} Q^{*}(\bm{s}_t,\bm{a}_t)
\end{align}
If the optimal Q-function is known, then the optimal action $\bm{a}^*_t$ can be extracted by choosing the action $\bm{a}_t$ that maximizes $Q^{{*}}(\bm{s}_t,\bm{a}_t)$ for state $\bm{s}_{t}$ (i.e. the optimal policy $\pi^*$ is retrieved),
\begin{equation}\label{eqn:Qmax}
	\bm{a}^*_t = \text{arg}\max_{\bm{a}_t \in \mathcal{A}} Q^{*}(\bm{s}_t,\bm{a}_t)
\end{equation}
without requiring knowledge of the environment dynamics. This last point is precisely why the Q-function enables model-free RL.

Lastly, the advantage function is defined to measure how advantageous the action is compared to the action drawn from policy,
\begin{align}
    A^{\pi}(\bm{s}_{t},\bm{a}_{t}) = Q^{\pi}(\bm{s}_{t},\bm{a}_{t}) - V^{\pi}(\bm{s}_{t})
\end{align}

\subsection{Actor-Critic}

Actor-critic is an approximate dynamic programming (ADP) method which solves dynamic programming heuristically. Importantly, the actor-critic approach allows for continuous state/action spaces by using a function approximator, e.g., a neural network. In RL, as well as in dynamic programming, the action is taken by a \textit{policy} to maximize the expected total discounted reward. By following a given policy and processing the rewards, one should estimate the expected return given states from the \textit{value function}. In the actor-critic approach, the \textit{actor} improves the policy based on the value function that is estimated by the \textit{critic}. We specifically focus on the policy gradient-based actor-critic algorithm in this work, and, in particular, the Proximal Policy Optimization (PPO) \cite{PPO} is considered. The critic is the parameterized value function $V_{\phi}$ and the actor is the parameterized policy $\pi_{\theta}$.

\subsubsection{Critic}
The role of the critic is to evaluate the current policy prescribed by the actor. The action is drawn from a Gaussian distribution, namely,
\begin{equation}\label{Gauss-policy}
    \bm{a}_{t} \sim \mathcal{N}(\mu, \sigma^2), \quad \text{where} \quad [\mu, \sigma] = f_{\text{DNN}}(\bm{s_t}; \theta),
\end{equation}
where a mean $\mu$ and standard deviation $\sigma$ computed from a deep neural network (DNN), $f_{\text{DNN}}(\bm{s_t}; \theta)  : \mathcal{S} \rightarrow \mathbb{R}^{2 }$, and the DNN is parameterized by weight vector $\theta$. After applying the action, we observe the reward $\bm{r}_{t+1}$ and the next state $\bm{s}_{t+1}$. For each time step $t$, the tuple $(\bm{s}_{t},\bm{a}_{t},r_{t+1},\bm{s}_{t+1})$ is stored in the buffer, $\mathcal{D}$. From a collected set of trajectories, the parameterized value function denoted by $V_{\phi}$ is updated to minimize the following loss function, $\mathcal{L}$:
\begin{align}
    \mathcal{L} = \frac{1}{N} \sum_{i \in \mathcal{D}} \sum_{t=0}^{N} \left( V_{\phi}\left(\bm{s}_{i,t}\right) - R_{i,t} \right)^{2}    
\end{align}
where $R_{i,t}$ is the total discounted reward at time $t$ in the i-th trajectory stored in the buffer. The critic network parameters, $\phi$, are updated numerically via gradient descent.

\begin{figure*}[t!]	
\includegraphics[width=\textwidth]{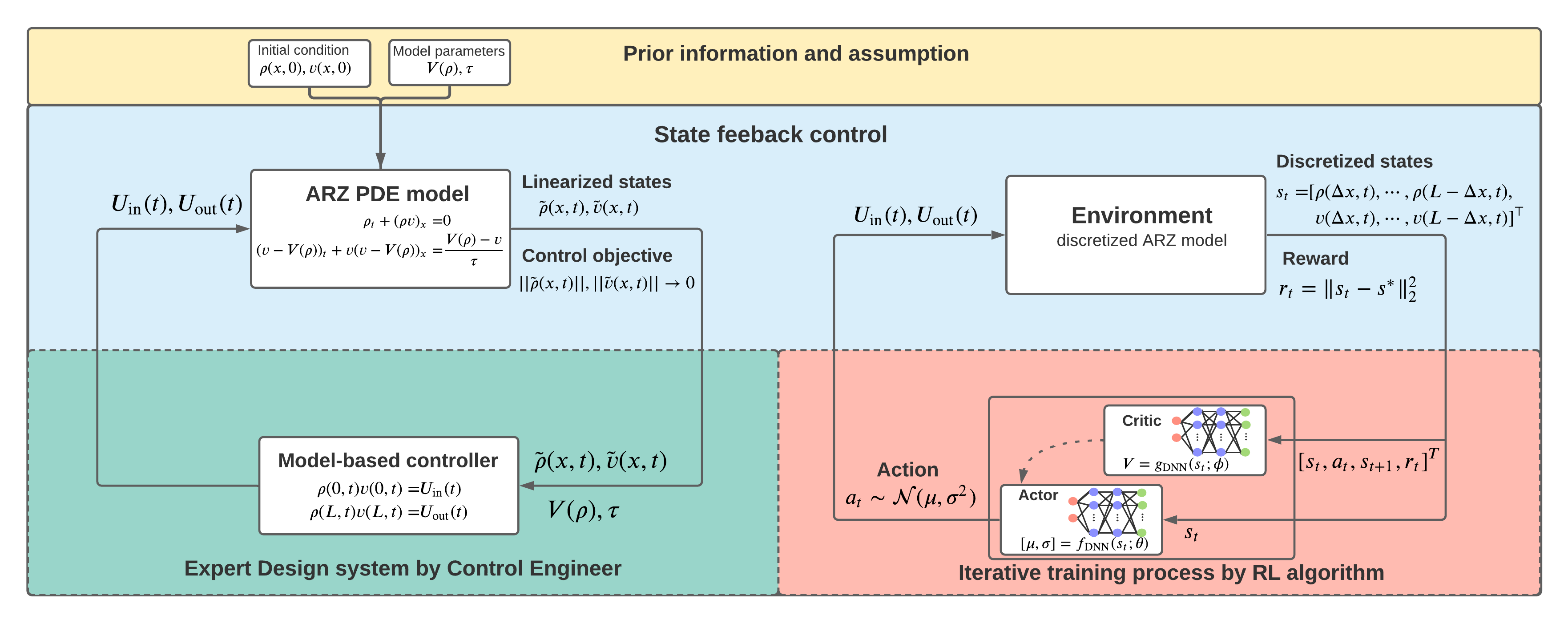}
	\caption{state-feedback control loop of ARZ PDE Lyapunov-based control scheme and model-free RL boundary control scheme}
	\label{fig:RL}
\end{figure*}

\subsubsection{Actor}
The actor is updated based on the value function estimates. The objective function for the actor is formulated in terms of the expected reward of policy $\pi_\theta$ and the advantage of $\pi_{\theta_{\text{old}}}$ \cite{Kakade,TRPO}:
\begin{align}
	\max_{\theta} \quad & \hat{\mathbb{E}}_{t} \left[ \frac{\pi_{\theta}(\bm{a}_{t}|\bm{s}_{t}) }{ \pi_{\theta_{\text{old}}}(\bm{a}_{t} | \bm{s}_{t}) } \hat{A}_{t} \right]
\end{align}
where the hats on $\hat{\mathbb{E}}_{t}$ signify a sample mean, and $\hat{A}_{t}$ indicates an estimated advantage function obtained from the critic. In \cite{TRPO}, the authors prove the expected reward corresponding to $\pi_\theta$ increases relative to $\pi_{\theta_{\text{old}}}$, if a distance measure between $\pi_\theta$ and $\pi_{\theta_{\text{old}}}$ is sufficiently bounded. This motivates the following trust region policy optimization algorithm:
\begin{align}
	\max_{\theta} \quad & \hat{\mathbb{E}}_{t} \left[ \frac{\pi_{\theta}(\bm{a}_{t}|\bm{s}_{t}) }{ \pi_{\theta_{\text{old}}}(\bm{a}_{t} | \bm{s}_{t}) } \hat{A}_{t} \right] \\
	\text{subject to} \quad &\hat{\mathbb{E}}_{t} \left[ \text{KL}[\pi_{\theta_{\text{old}}}(\cdot | \bm{s}_{t}),\pi_{\theta}(\cdot | \bm{s}_{t})] \right] \leq \delta,
\end{align}
where $\theta_{\text{old}}$ is the vector of policy parameters before the update. KL-divergence measures the difference between the old policy and current policy. The constraint ensures that the new policy does not deviate from the old policy by $\delta$. 

In this work, we adopt the PPO reinforcement learning algorithm \cite{PPO}, which is based on trust region policy optimization (TRPO) \cite{TRPO}. The PPO algorithm similarly limits the new policy from being excessively far from the previous one. However, it does so with a modified objective that penalizes changes to the policy that move $r_{t}(\theta) = \pi_{\theta}(\bm{a}_{t}|\bm{s}_{t}) / \pi_{\theta_{\text{old}}}(\bm{a}_{t}|\bm{s}_{t})$ away from 1. The key idea is to use probability clipping, as follows:
\begin{align}
	\max_{\theta} \quad & \hat{\mathbb{E}}_{t} \big[  \min(r_{t}(\theta)\hat{A}_{t}, \text{clip}(r_{t}(\theta), 1-\varepsilon, 1+\varepsilon ) \hat{A}_{t} \big]. \label{eqn:actor_update}
\end{align}
The main idea of PPO is to modify the objective by clipping the probability ratio. This removes the incentive for moving $r_{t}$ outside of the interval $[1-\varepsilon, 1+\varepsilon]$. With this clipping method, the lower bound of objective function is maximized. Readers are referred to \cite{PPO} for more details.

\subsubsection{Proximal Policy Optimization}

We adopt a policy gradient-based approach to obtain a continuous-valued stochastic control policy. Mathematically, the goal is to find:
\begin{equation}
    \theta^\star = \arg \max_\theta \mathbb{E} \left[ \mathcal{R}_{t} \right] \label{eqn:RL_Obj}
\end{equation}
where the expectation is taken w.r.t. $\mathcal{P}(\bm{s}_{t+1} | \bm{s}_{t}, \bm{a}_{t})$ and $\pi_\theta( \bm{a}_t | \bm{s}_t)$, and $\theta$ parameterizes the control policy distribution. Policy gradient methods essentially solve \eqref{eqn:RL_Obj} via gradient ascent. The key challenge is estimating the gradient since it is computationally intractable to compute it exactly.

One can re-formulate this optimization problem \eqref{eqn:RL_Obj} in terms of the expected reward of policy $\pi_\theta$ and the advantage of $\pi_{\theta_{\text{old}}}$ \cite{Kakade,TRPO}:
\begin{align}
	\max_{\theta} \quad & \hat{\mathbb{E}}_{t} \left[ \frac{\pi_{\theta}(\bm{a}_{t}|\bm{s}_{t}) }{ \pi_{\theta_{\text{old}}}(\bm{a}_{t} | \bm{s}_{t}) } \hat{A}_{t} \right]
\end{align}
where the hats on $\hat{\mathbb{E}}_{t}$ signify a sample mean, and $\hat{A}_{t}$ indicates an estimated advantage function from simulations. In \cite{TRPO}, the authors prove the expected reward corresponding to $\pi_\theta$ increases relative to $\pi_{\theta_{\text{old}}}$, if a distance measure between $\pi_\theta$ and $\pi_{\theta_{\text{old}}}$ is sufficiently bounded. This motivates the following trust region policy optimization algorithm:
\begin{align}
	\max_{\theta} \quad & \hat{\mathbb{E}}_{t} \left[ \frac{\pi_{\theta}(\bm{a}_{t}|\bm{s}_{t}) }{ \pi_{\theta_{\text{old}}}(\bm{a}_{t} | \bm{s}_{t}) } \hat{A}_{t} \right] \\
	\text{subject to} \quad &\hat{\mathbb{E}}_{t} \left[ \text{KL}[\pi_{\theta_{\text{old}}}(\cdot | \bm{s}_{t}),\pi_{\theta}(\cdot | \bm{s}_{t})] \right] \leq \delta,
\end{align}
where $\theta_{\text{old}}$ is the vector of policy parameters before the update. KL-divergence measures the difference between the old policy and current policy. The constraint ensures that the new policy does not deviate from the old policy by $\delta$. 

Both the RL and Lyapunov-based approaches provide state-feedback controllers for the ARZ PDE model. Figure \ref{fig:RL} compares the RL and Lyapunov-based approaches as signal flow diagrams. The Expert Design system developed by a Control Engineer is replaced with an iterative learning process by a RL algorithm. In addition, the prior assumption of the initial condition and model parameters are not required in the RL approach. In the next section, we will conduct several numerical simulations to compare these two methodologies.

\section{Comparative study and simulation}
In this section, we numerically test the control designs of the Lyapunov-based controllers and the RL controller. Two settings are considered: (i) full knowledge of the system dynamics, and (ii) partial knowledge of the system. In the full knowledge setting, both the Lyapunov-based and model-free RL controllers have perfect knowledge of the model with known model parameters and steady state conditions. For the partial knowledge setting, both controllers have partial knowledge of the model, i.e., the true steady state has deviated from the ones used in the model. For each setting, we test the model-free RL controller with the following cases: (i) outlet control, (ii) inlet control, and (iii) outlet \& inlet control. For each case, one Lyapunov-based controller is used to evaluate the performance of the RL controller.

\subsection{Comparative study with full knowledge of system dynamics} \label{det}
We simulate traffic on a freeway segment with length $L = 500 \, \rm{m}$ for a time period $T = 240\,  \rm{s} = 4  \, \rm{min}$.
 Steady states density $\rho^\star= 120 \, \rm{veh}/\rm{km}, v^\star = 36\, \rm{km}/\rm{h}$ are chosen given the maximum density $\rho_{ m} = 160 \, \rm{veh}/\rm{km}$, and maximum velocity $v_{ m} =40\,  \rm{m}/\rm{s}$. The steady state traffic is lightly congested such that $\lambda_1 = 10 \, \rm{m/s},~ \lambda_2 = -20\,  \rm{m/s}$, which satisfy the conditions in \eqref{char1} -- \eqref{char2}. This configuration represents heterogeneous propagation of oscillatory waves of density and velocity in congested traffic.
We assume sinusoidal initial conditions:
\begin{align}
\rho(x,0) &=  0.1 \sin\left(\frac{3\pi x}{L}\right)\rho^\star + \rho^\star, \label{ic1} \\
v(x,0) &= -0.1 \sin\left(\frac{3\pi x}{L}\right) v^\star + v^\star. \label{ic2}
\end{align}

\begin{figure*}[ht!]
    \centering
        \begin{subfigure}[h]{0.485\textwidth}
        \centering
        \includegraphics[trim = 0.0mm 0.0mm 0.0mm 20.0mm, clip, width=0.9\textwidth]{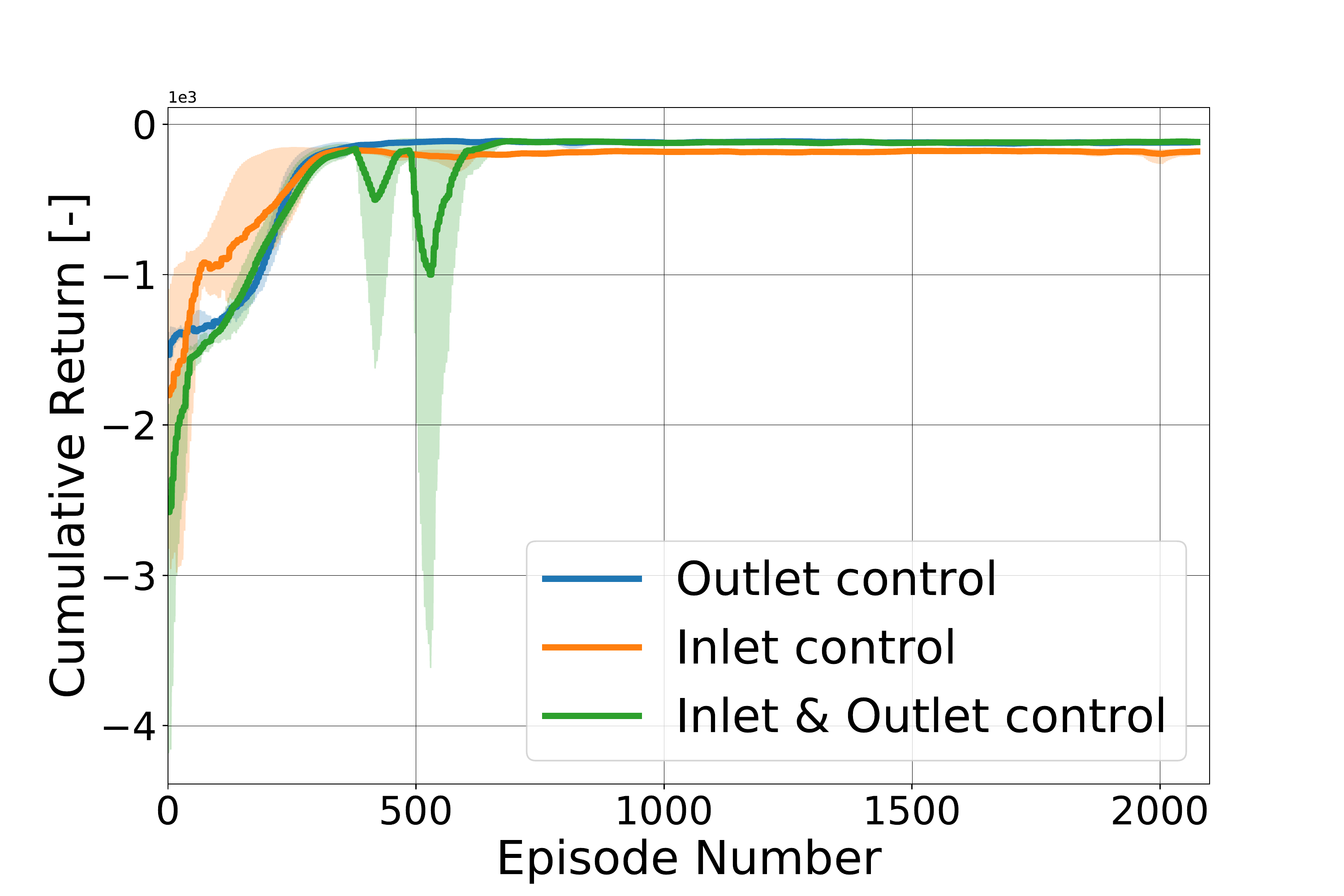}
        \caption{The learning curve of RL controllers with different initializations of the actor-critic network.}
        \vspace{-1ex}
        \label{fig:det_learning_curve}
    \end{subfigure}
    \hfill
    \begin{subfigure}[h]{0.49\textwidth}
        \centering
        \includegraphics[trim = 10.0mm 0.0mm 0.0mm 10.0mm, clip, width=0.9\textwidth]{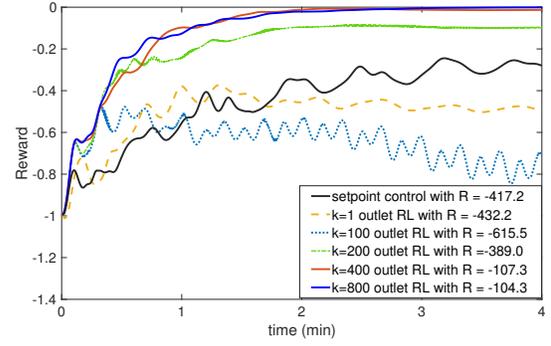}
        \caption{The reward evolution for the closed-loop system with outlet RL controller that is obtained for different episode numbers k.}
        \vspace{0ex}
        \label{fig:epi_reward}
    \end{subfigure}%
      \caption{The learning curves for RL and reward evolution in the learning process.}
    \vspace{0ex}
    \label{fig:det_reward_curves}
\end{figure*}
\subsubsection{Learning process of RL controllers} The learning process for RL controllers is illustrated in Fig.~\ref{fig:det_reward_curves}. The learning curve for RL corresponding to each control scheme is presented in Fig.~\ref{fig:det_learning_curve}. The evolution of the cumulative reward $R_t$, defined in \eqref{eqn:reward_fct}, reflects the overall learning performance. In addition, several episodes of the simulation in the learning process are compared more closely in Fig.~\ref{fig:epi_reward}. We plot the reward values $r_t$, defined in \eqref{reward}, over several hundreds of iterations of testing candidate RL outlet feedback laws on ARZ model simulations over a time window of 4 minutes for each simulation.

In Fig.~\ref{fig:det_learning_curve}, the shaded area represents min/max of different actor-critic parameterizations while the bold line measures the average of the performance during training. Notice that the training process converges over different numbers of episodes across the controllers. This motivates us to set a large enough number of episodes in order to make sure that RL converges. We obtain the actor-critic parameters by the end of training. It is also observed that there is a spike in cumulative reward for the inlet and outlet RL controller at the episode around 500 which indicates that the closed-loop system at the iteration becomes unstable and states diverge. For example, if density state is greater than the maximum density, there will be collision in the ARZ model.  Therefore, this training process cannot perform in a collision free fashion if applied in real traffic. 

In Fig.~\ref{fig:epi_reward}, we choose the outlet control case to illustrate the RL learning process. We plot the reward evolution for the closed-loop system with RL outlet controllers that are obtained from training episodes 1, 100, 200, 400 and 800 respectively. In addition, the reward of the setpoint control is also plotted as a baseline comparison. It turns out that the initial reward is worse than the baseline performance.  After 100 iterations of training, the outlet RL controller plotted with dotted blue line performs even worse with reward decreasing with time and the lowest cumulative reward R = -615.5. However, the RL performance improves after the iteration number increases from 200 to 400 and reaches its best at episode 800. As demonstrated in this figure, the learning process does not guarantee a monotonically improving performance. Certain episodes can be even worse than the baseline case due to some iterations of unsuccessful training.

\begin{figure*}[ht!]
\begin{subfigure}[a]{.6\textwidth}
	\centering
	\includegraphics[width=1\textwidth]{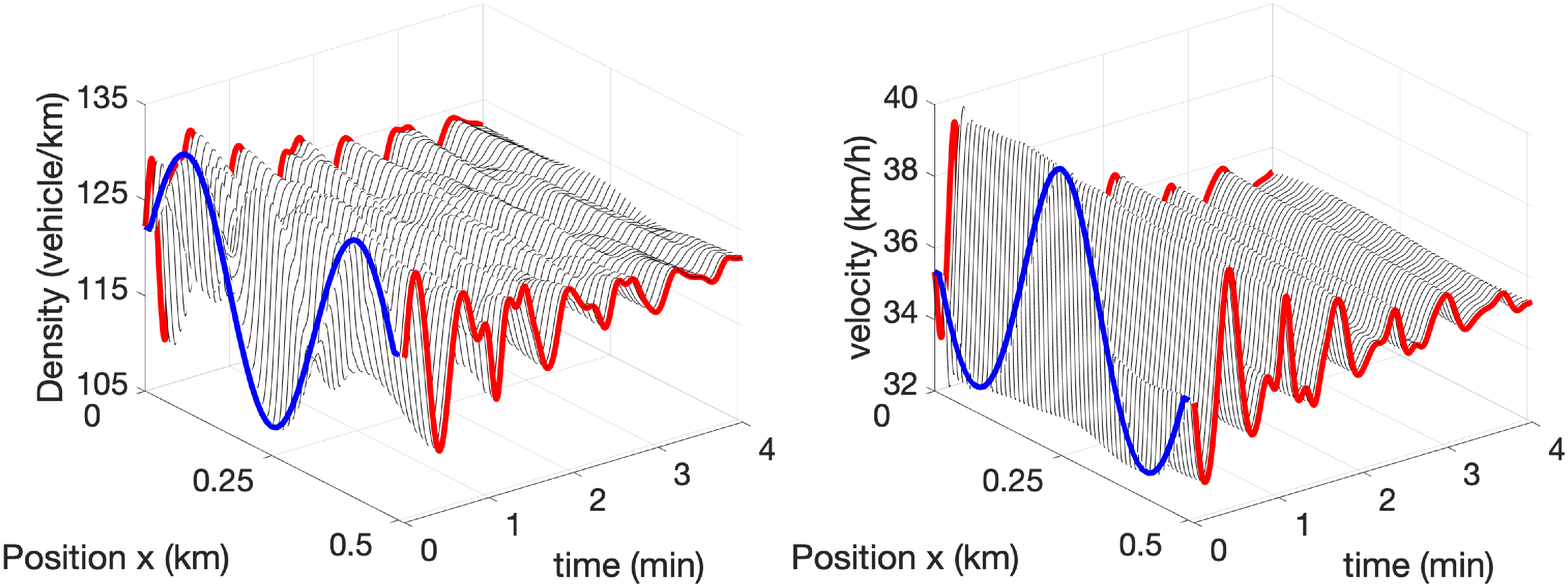}
	\caption{Density and velocity of ARZ model under setpoint control.}
	\label{fig:openloop_dens_velo}
\end{subfigure}
\begin{subfigure}[a]{.4\textwidth}
	\centering
	\includegraphics[width=6cm]{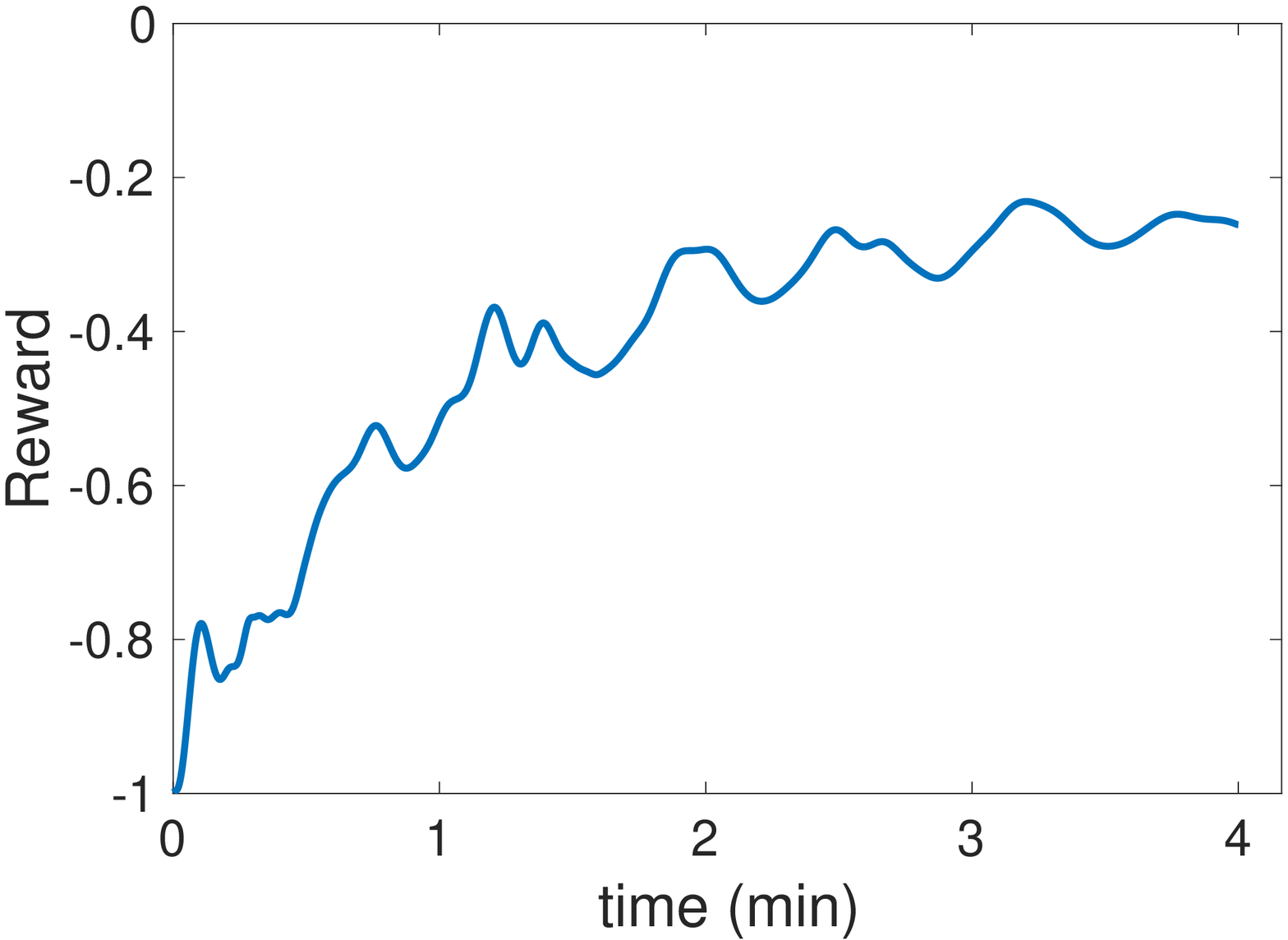}
	\label{fig:openloop_reward}
	  \caption{Reward evolution for setpoint control}
\end{subfigure}
\caption{Density and velocity evolution in space and time with setpoint inputs.}
\label{fig:openloop}
\end{figure*}

\begin{figure*}[ht!]
\hspace*{-0.2cm}
\begin{subfigure}[a]{.6\textwidth}
  \centering
  \includegraphics[width=1.05\linewidth]{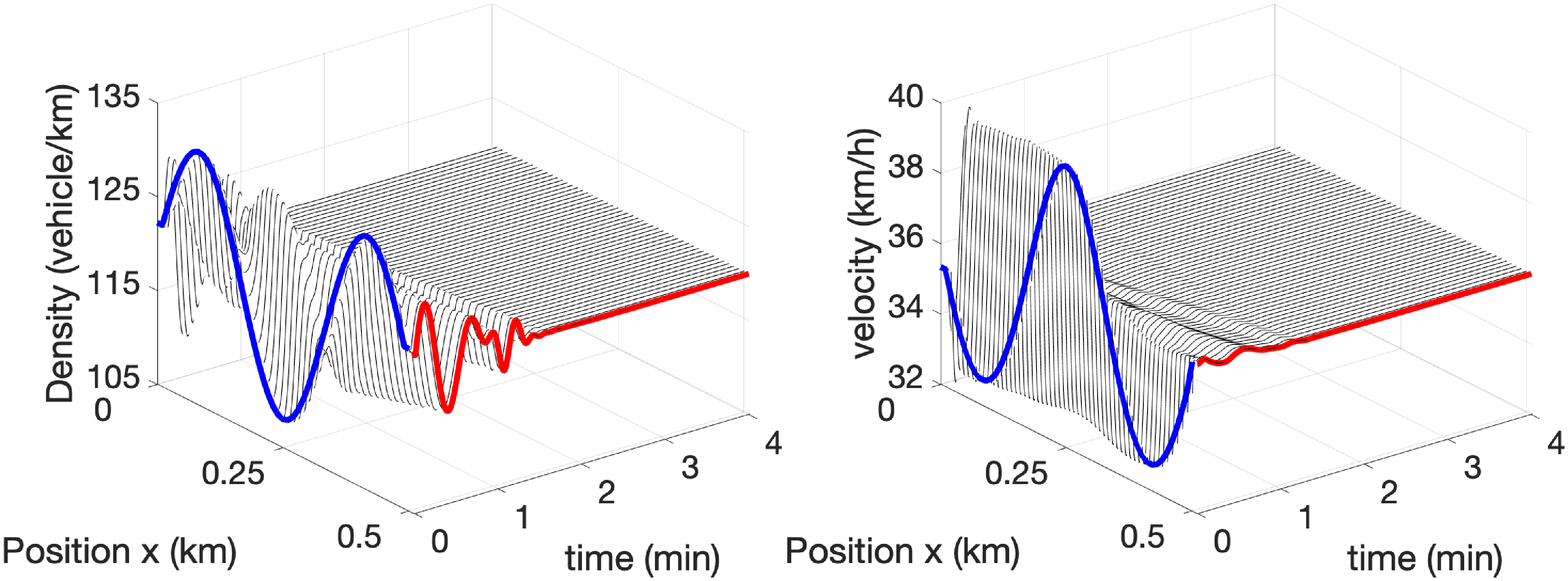}  
  \caption{closed-loop with {\bf outlet backstepping controller}}
  \label{fig:bkst}
 \hspace*{-0.15cm}
  \includegraphics[width=1.05\linewidth]{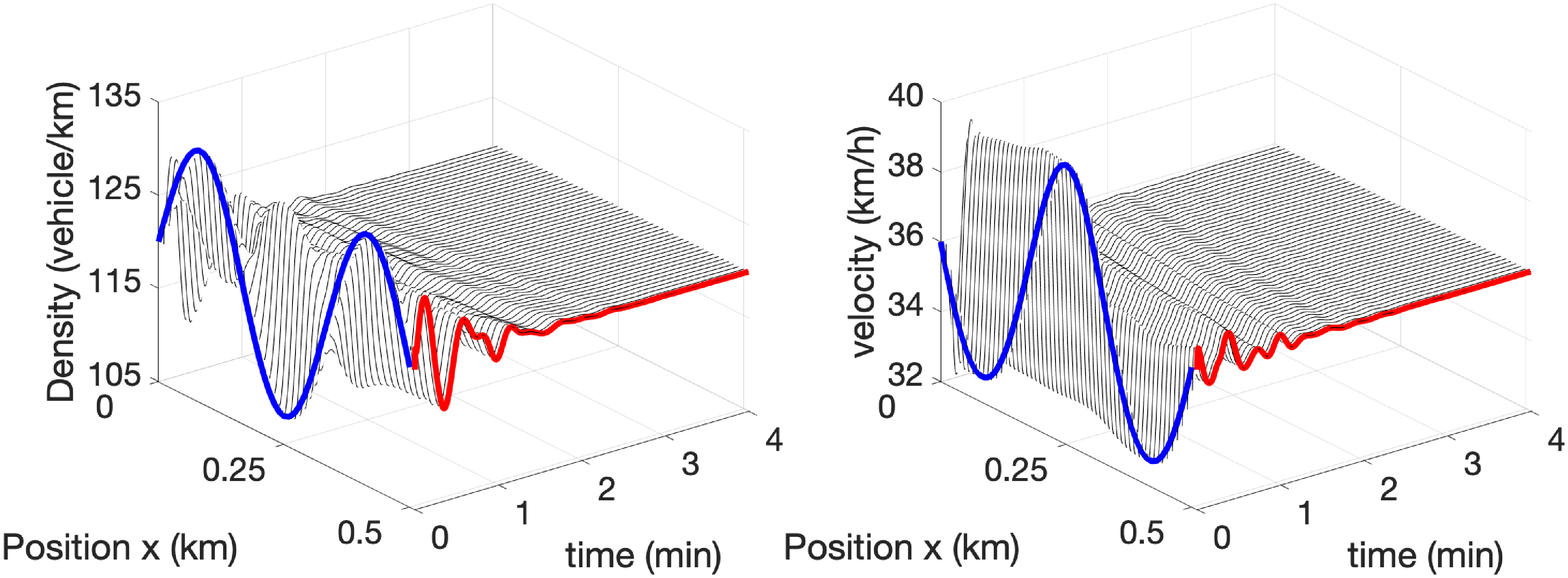}  
  \caption{closed-loop with {\bf outlet RL controller}}
  \label{fig:RL_bkst}
\end{subfigure}
\hspace*{-0.3cm}
\begin{subfigure}[c]{0.4\textwidth}  
\vspace{0.1cm}
       \centering 
            \includegraphics[width=0.9\textwidth]{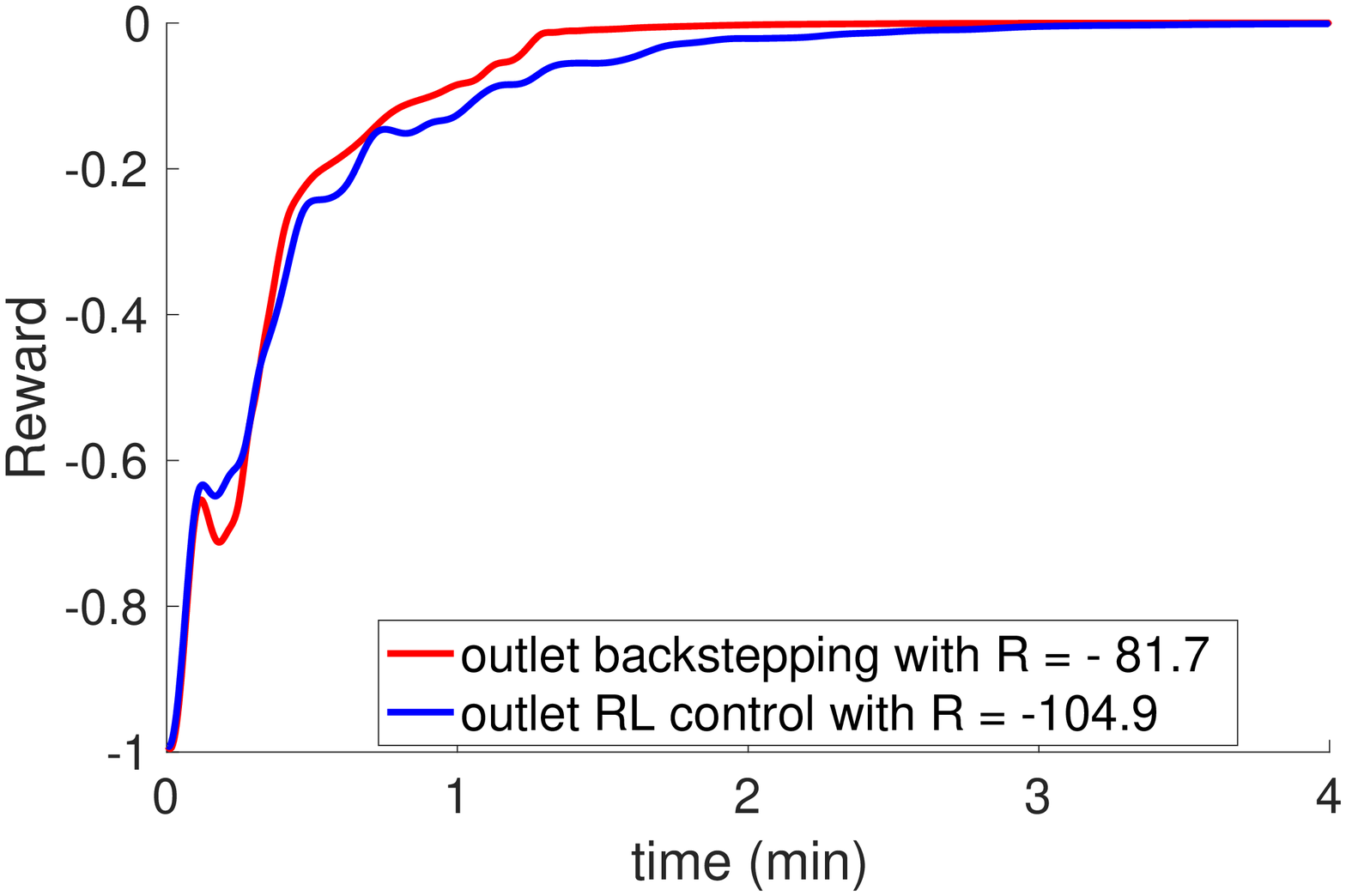}
            \caption{Reward evolution for closed-loop with backstepping and RL controller at outlet }
            \label{fig:reward_bkst}
\vspace{0.3cm}
            \centering 
            \includegraphics[width=0.9\textwidth]{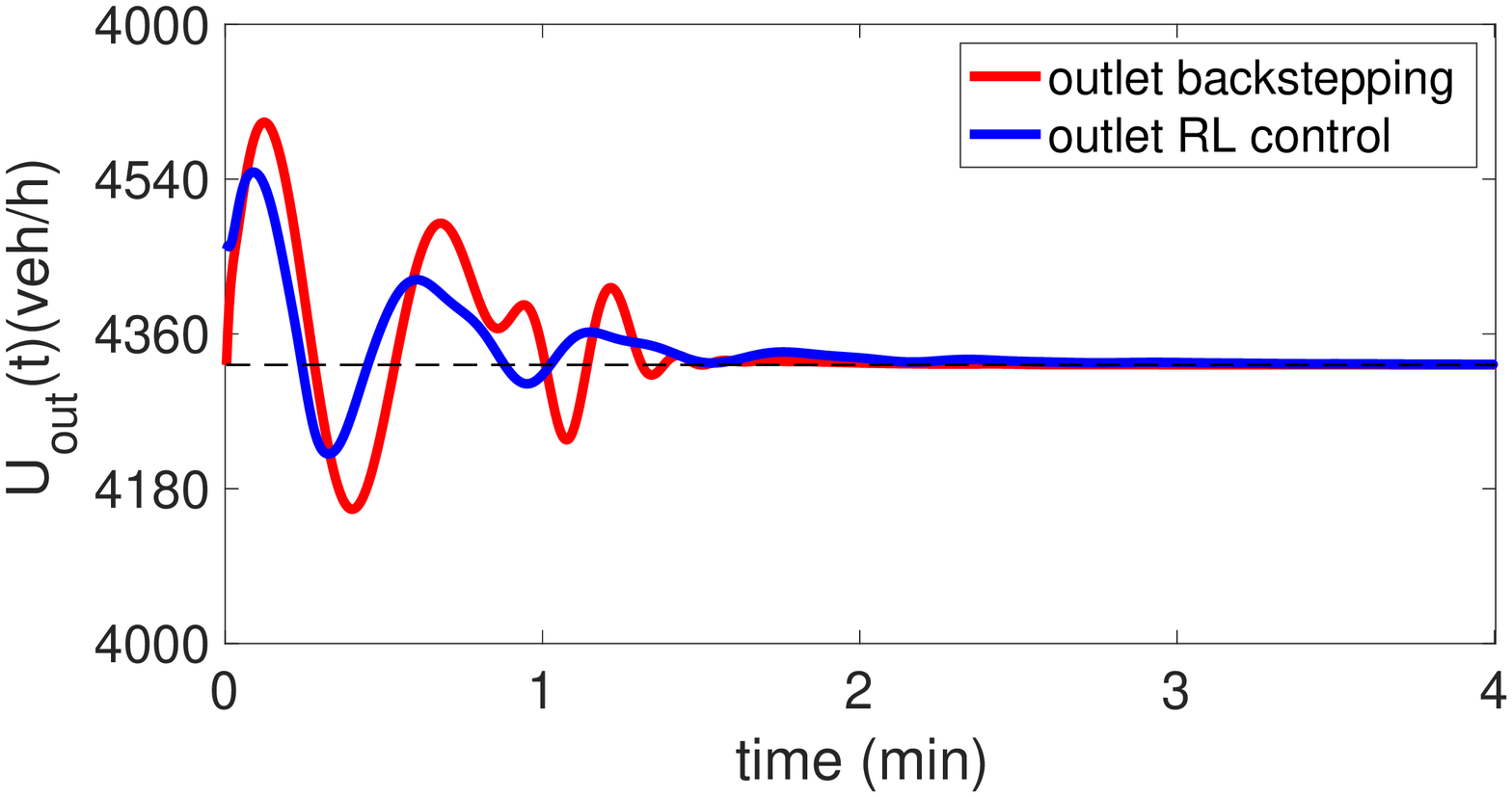}
            \caption{Backstepping and RL control inputs at outlet} 
            \label{fig:control_bkst}
        \end{subfigure}
\caption{Density and velocity evolution in space and time with control input highlighted with red at outlet }
\label{fig:state_bkst}
\end{figure*}

\begin{figure*}[ht!]
\hspace*{-0.2cm}
\begin{subfigure}[a]{.6\textwidth}
  \centering
  \includegraphics[width=1.05\linewidth]{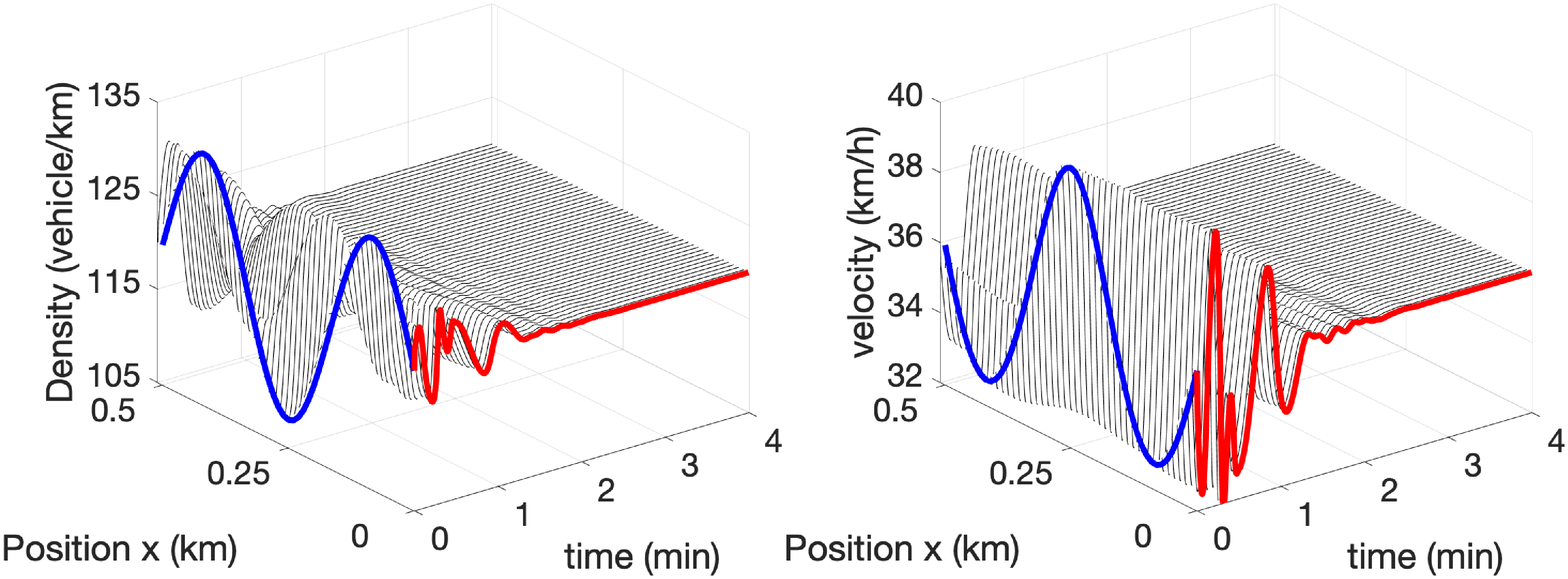}  
 \caption{closed-loop with {\bf inlet P controller}}
  \label{fig:P}
 \hspace*{-0.15cm}
  \includegraphics[width=1.05\linewidth]{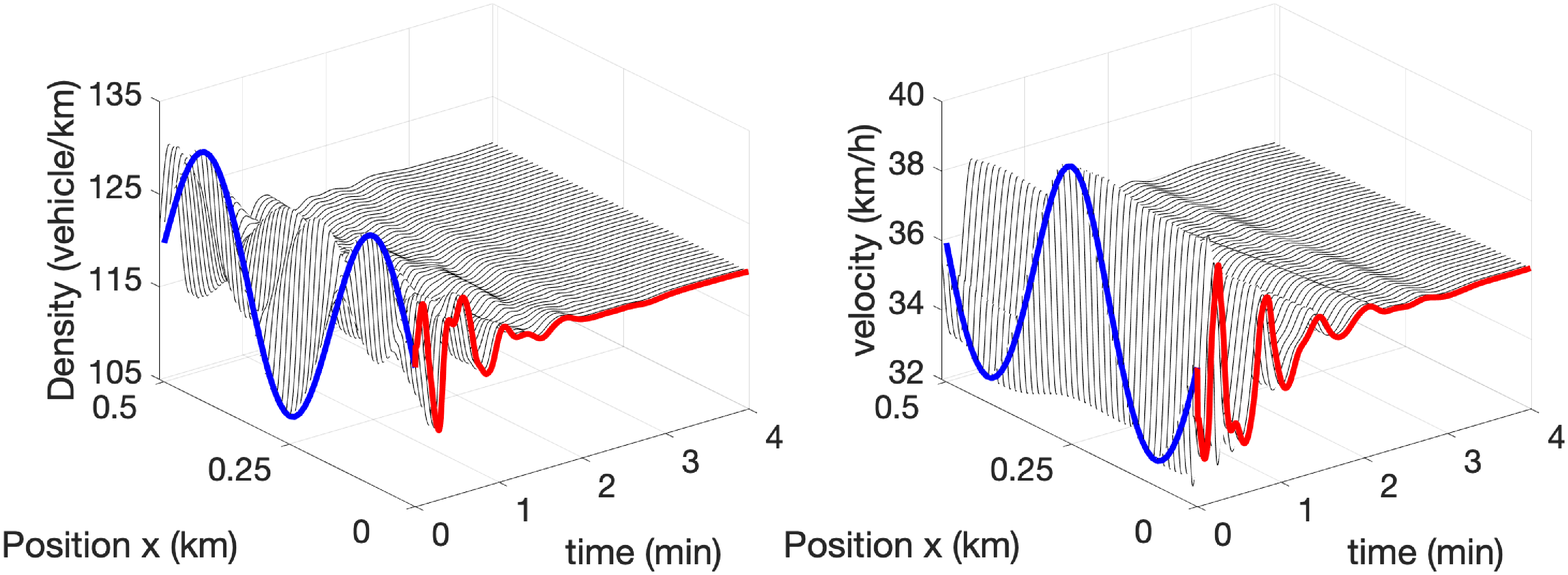}
 \caption{closed-loop with {\bf inlet RL controller}}
  \label{fig:RL_p}
\end{subfigure}
\hspace*{-0.3cm}
\begin{subfigure}[c]{0.4\textwidth}  
\vspace{0.1cm}
       \centering 
            \includegraphics[width=0.85\textwidth]{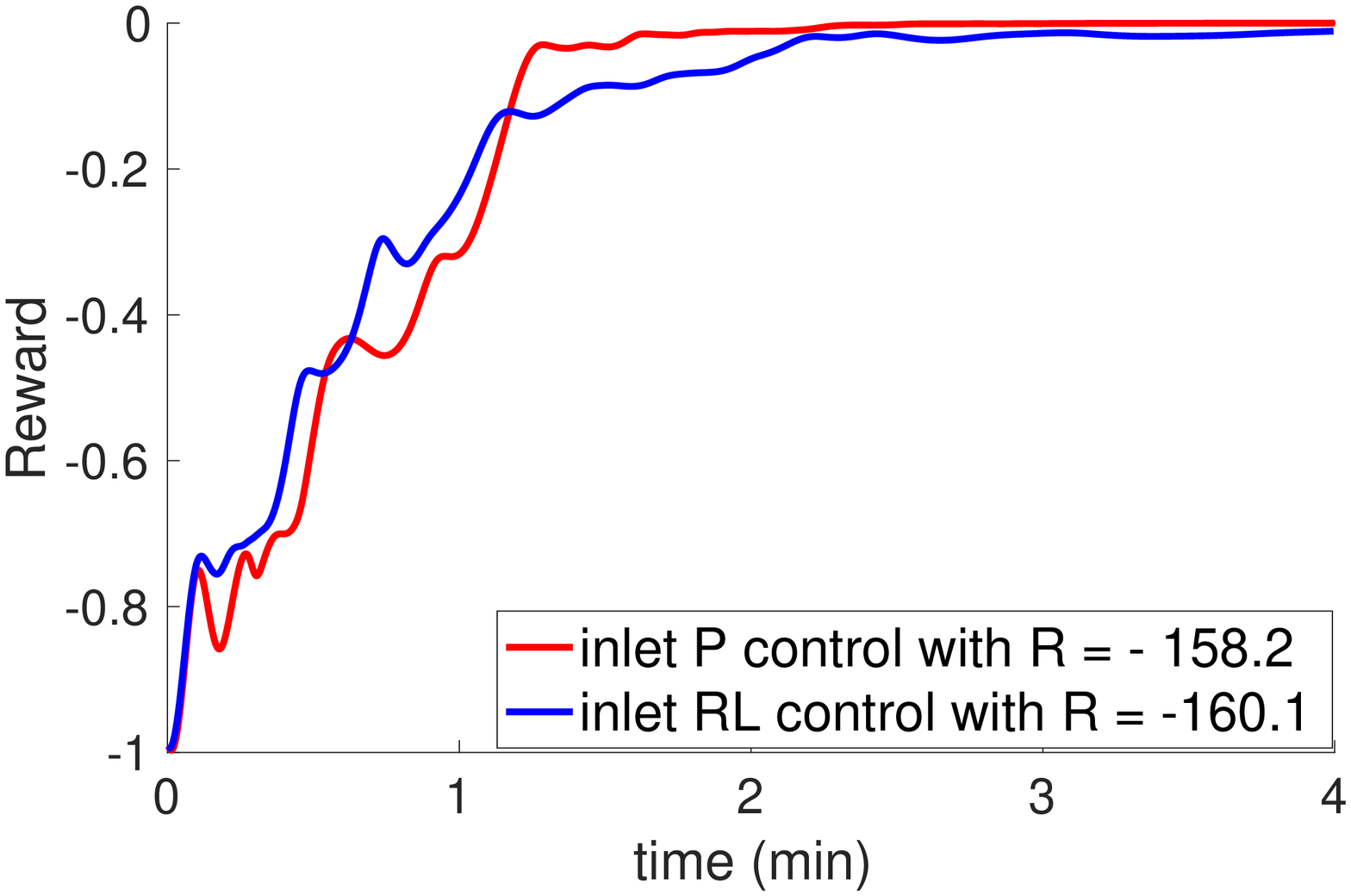}
            \caption{Reward evolution for closed-loop with P and RL controller at inlet.}
            \label{fig:reward_p}
\vspace{0.3cm}
            \centering 
            \includegraphics[width=0.9\textwidth]{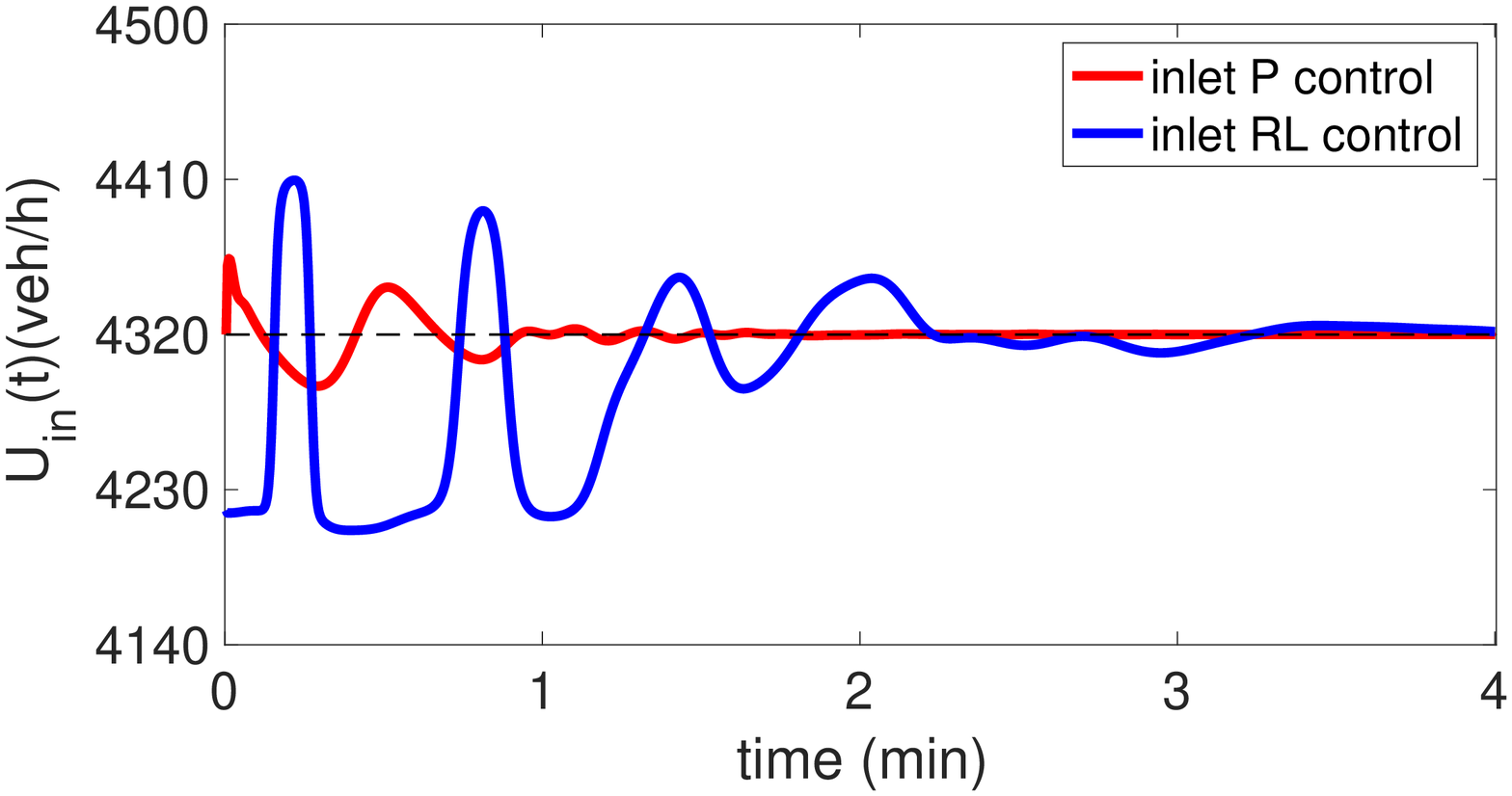}
            \caption{ P and RL control inputs at inlet.} 
            \label{fig:control_p}
        \end{subfigure}
\caption{Density and velocity evolution in space and time with control input highlighted with red at inlet.}
\label{fig:state_p}
\end{figure*}

\begin{figure*}[ht!]
\hspace*{-0.2cm}
\begin{subfigure}[a]{.6\textwidth}
  \centering
  \includegraphics[width=1.05\linewidth]{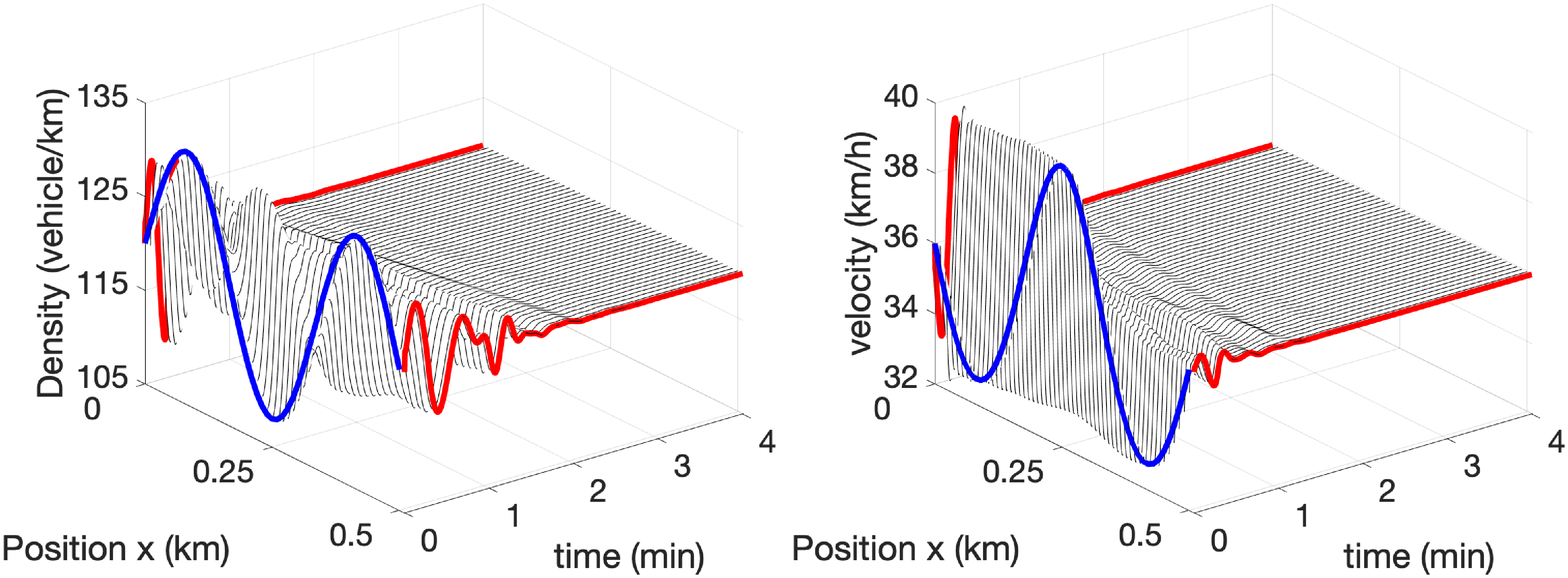}  
 \caption{closed-loop with {\bf PI controllers} at inlet and outlet}
  \label{fig:PI}
 \hspace*{-0.15cm}
  \includegraphics[width=1.05\linewidth]{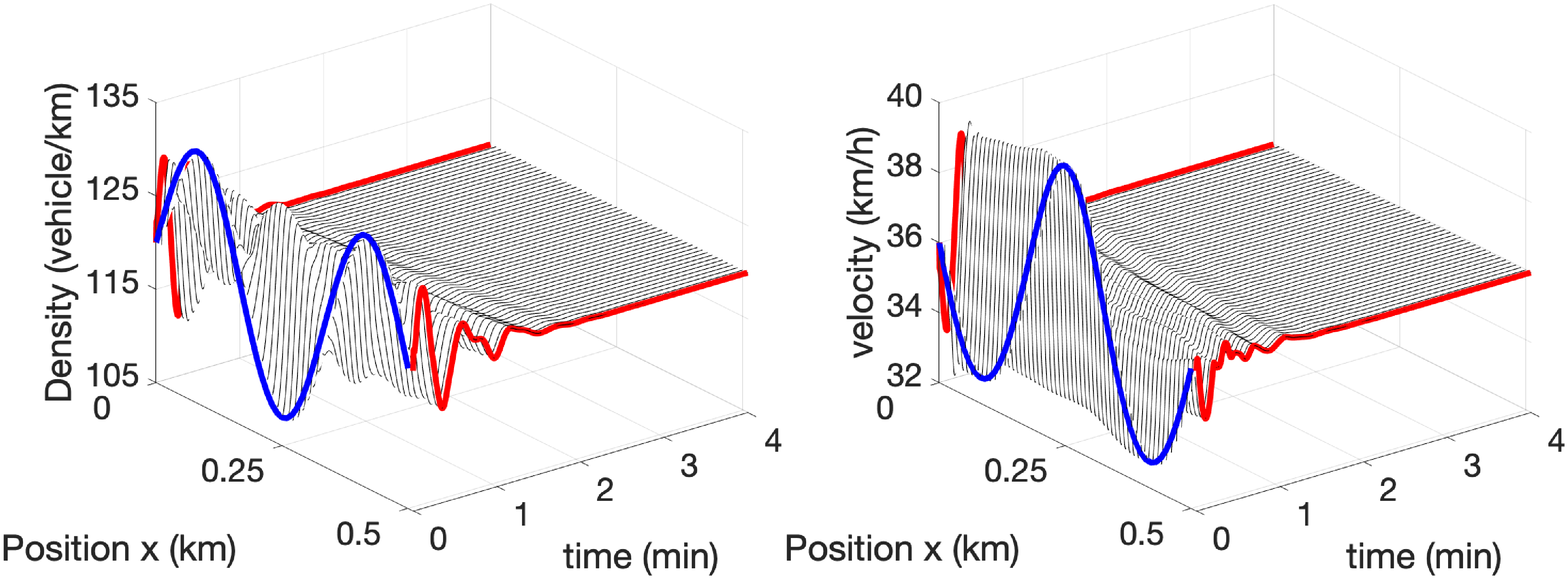}
  \caption{closed-loop with {\bf RL controllers} at inlet and outlet}
  \label{fig:RL_PI}
\end{subfigure}
\hspace*{-0.3cm}
\begin{subfigure}[c]{0.4\textwidth}  
\vspace{0.1cm}
       \centering 
            \includegraphics[width=0.85\textwidth]{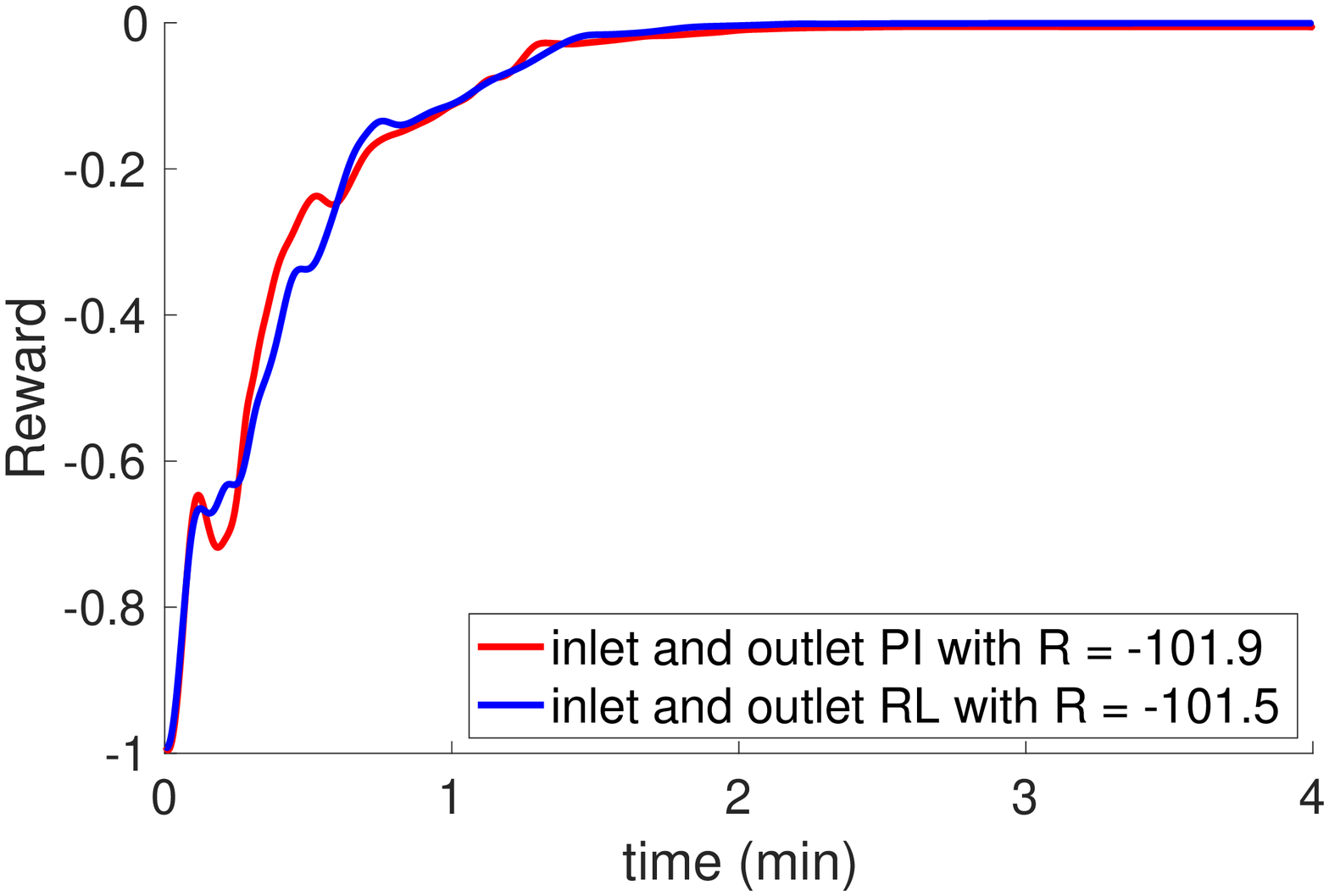}
            \caption{Reward evolution for closed-loop with PI and RL controller at inlet and outlet }
            \label{fig:reward_PI}
\vspace{0.1cm}
            \centering 
            \includegraphics[width=0.8\textwidth]{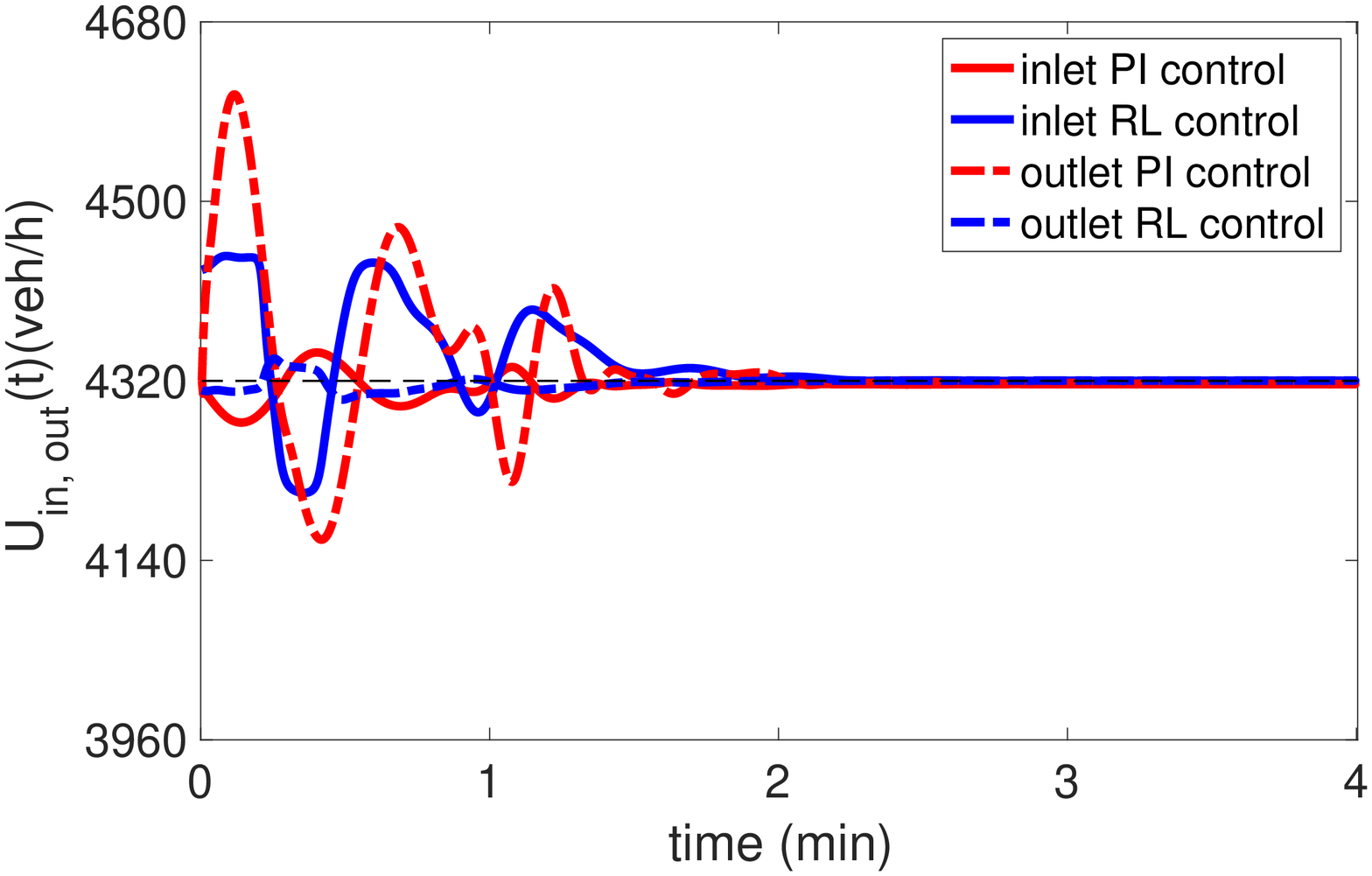}
            \caption{Inlet and outlet PI and RL control inputs} 
            \label{fig:control_PI}
        \end{subfigure}
\caption{Density and velocity evolution in space and time with {\bf inlet and outlet control input highlighted with red}}
\label{fig:state_PI}
\end{figure*}

\subsubsection{State evolution, reward and control inputs}
In Figs.~\ref{fig:openloop}--\ref{fig:state_PI}, the initial condition is highlighted in blue. The actuated inlet and/or outlet boundaries are highlighted in red. In other words, the red curves visualize the control inputs for the setpoint control, backstepping control, P control, PI control and the RL controllers in each case. Specifically, the closed loop results are compared between the outlet backstepping controller and the outlet RL controller, as shown in Fig.~\ref{fig:state_bkst}, between inlet P controller and inlet RL controller as shown in Fig.~\ref{fig:state_p}, between inlet and outlet PI controller and inlet and outlet RL controller as shown in Fig.~\ref{fig:state_PI}. 

Figure~\ref{fig:openloop} plots the evolution of density and velocity under setpoint control. We observe the persisting oscillations, although they appear lightly damped. We regard the setpoint control as a baseline case against which the subsequent feedback control designs should outperform. In setpoint control, the outlet and inlet boundary flow rates are controlled to be constant \eqref{bc1},\eqref{bc2}. As highlighted with red in~\ref{fig:openloop}, the product of density and velocity values is constant at inlet and outlet boundary whereas the boundary states of density and velocity oscillate over time.

Fig.~\ref{fig:state_bkst} presents the closed-loop result of the outlet backstepping controller and outlet RL controller. In Fig.~\ref{fig:bkst} and Fig.~\ref{fig:RL_bkst}, the states are stabilized to spatially uniform steady state values by the outlet backstepping controller and outlet RL controller, respectively.
An interesting finding from comparing between the two controllers is that RL learns a policy which produces a control input (red line in Fig.~\ref{fig:RL_bkst}) that closely replicates the backstepping control input (red line in Fig.~\ref{fig:bkst}). The RL policy is developed without explicit knowledge of the differential equations and parameters. It is trained iteratively on the nonlinear simulation model. As a state-feedback controller, the RL policy produces an action for boundary input given the current state. 
In contrast, the PDE backstepping state-feedback control law is obtained by rigorous theoretical control design assuming perfect knowledge of the model. As shown in Fig.~\ref{fig:control_bkst}, both methods yield similar control input trajectories. However, RL underperforms relative to backstepping in terms of the instantaneous reward over time. The convergence of the $L^2$ spatial norm represents the stabilization of the closed-loop system. Unlike the finite-time convergence for the backstepping controller shown with red line, the RL controller's reward in blue takes a longer time -- about $3 \, \rm{min}$ for convergence. The cumulative reward for the RL outlet controller is $R = -104.9$ while the cumulative reward of backstepping is $R = -81.7$.

In the same fashion, we compare closed-loop results for inlet P controller and inlet RL controller in Fig.~\ref{fig:state_p}. The incoming flow rate is actuated either with the P controller or RL controller. The control inputs at the inlet are highlighted with red in Fig.~\ref{fig:P} and Fig.~\ref{fig:RL_p}. As shown in the Fig.~\ref{fig:control_p}, the inlet RL controller needs more than $4 \,\rm min$ to converge, longer than the P control. Although the inlet RL controller almost recovers the performance of the inlet P controller, we can see from Fig.~\ref{fig:reward_p} that the RL control input is quite different than the P control input, requiring larger control effort and longer convergence time.

Figure~\ref{fig:state_PI} describes the closed-loop results for PI controller and RL controller, actuating the incoming and outgoing flow rates. We observe that the RL controller achieves a better reward evolution in Fig.~\ref{fig:reward_PI}. The convergence of the closed-loop system with RL is faster than that of the PI controller, and with a slightly smaller cumulative reward. It could also be noted that control efforts are smaller for the RL controller in terms of magnitude, as shown in Fig.~\ref{fig:reward_PI}. These results will vary, however, with the PI controller gain selection.

\begin{figure}[t!]
  \includegraphics[width=1.05\linewidth]{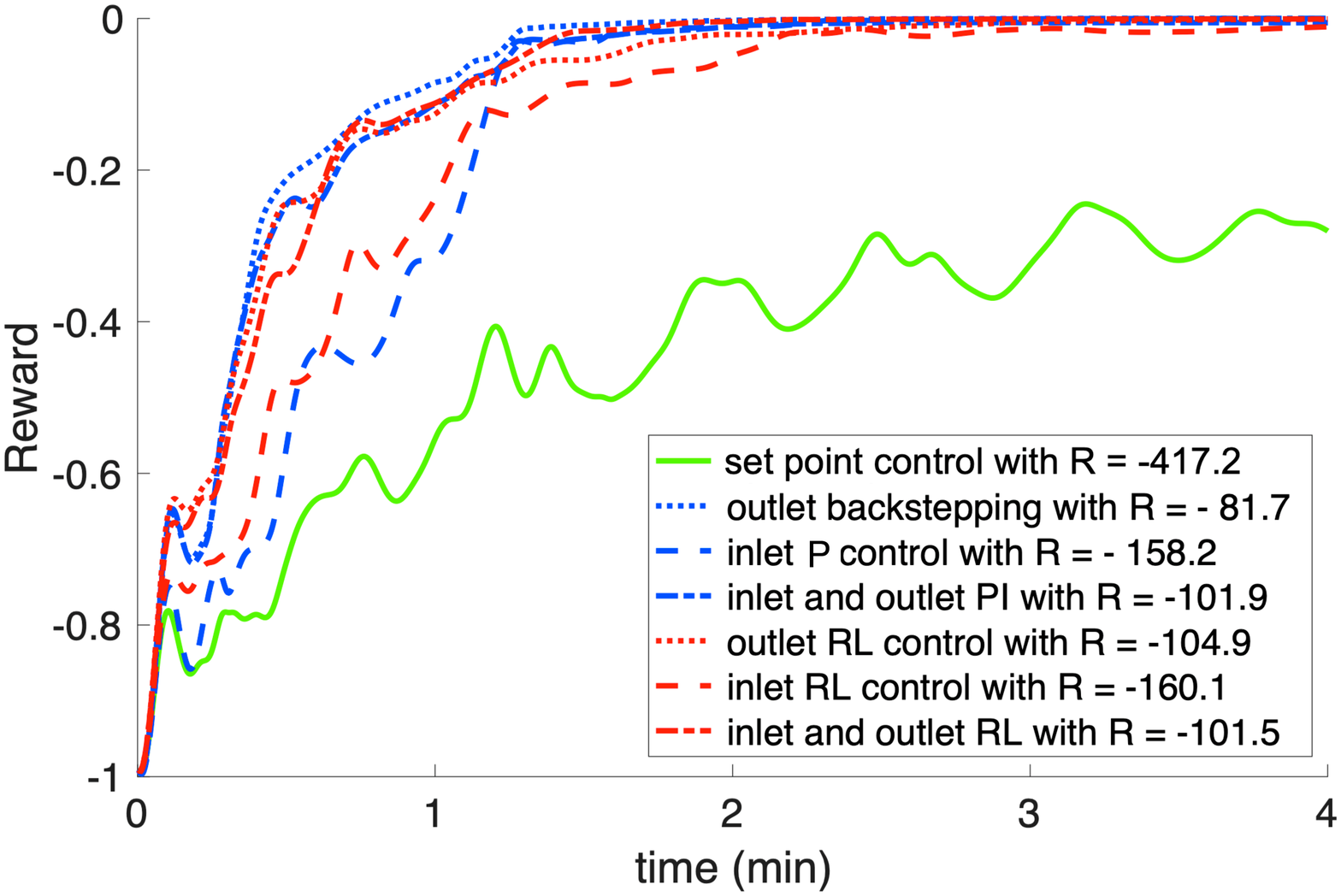}
            \caption{The reward evolution of the Lyapunov-based approaches (blue), the corresponding RL controllers after 2000 episodes of iterative training (red).} 
            \label{fig:reward_comparison}
\end{figure}

Figure~\ref{fig:reward_comparison} includes all the closed-loop results with the Lyapunov-based controllers and their corresponding RL controllers. All of the Lyapunov-based feedback controllers and RL controllers outperform the setpoint controllers. The RL agent achieves sub-optimal performance compared to Lyapunov-based approaches. In addition, backstepping outperforms the P and PI controllers for stabilization of the traffic flow. RL controllers almost recover the closed-loop results of the Lyapunov-based design in the outlet control or inlet control cases. RL outperforms the PI controllers when both the inlet and outlet boundaries are controlled. The aforementioned results are uncertain model parameters, which we will discuss in the next subsection. Note that the reward in \eqref{reward} is defined as the $L^2$ spatial norm of the state values without considering control effort. The RL controllers were design for comparison with the stabilizing Lyapunov-based controllers, and thus control effort was not penalized. However, RL controllers can be trained to consider control effort by modifying \eqref{reward} with a control effort penalty.

\subsubsection{Other performance measures}
Common traffic performance indices include total travel time (TTT), fuel consumption, and travel comfort. These indices are calculated for the closed-loop system with either the Lyapunov-based controllers or the RL controllers. The performance indices are given by
	\begin{align}
	J_{\rm TTT} &= \int_0^{T}\int_0^L\rho(x,t) d x d t, \\
	J_{\rm fuel} &= \int_0^{T}\int_0^L\max\{0,b_0+b_1v(x,t)+b_3v(x,t) \notag \\
	&+b_4v(x,t)a(x,t)\}\rho(x,t) d x d t, \\
	J_{\rm comfort} &= \int_0^{T}\int_0^L(a(x,t)^2+a_t(x,t)^2)\rho(x,t) d x d t
	\end{align}
where $a(x,t)$ is defined as the local acceleration $a(x,t)=v_t(x,t)+v(x,t)v_x(x,t)$ and $b_i$ are constant coefficients chosen as $b_0=25\cdot 10^{-3}{[\rm l/s]}$, $b_1=24.5\cdot10^{-6}{[\rm{l/m}]}$, $b_3=32.5\cdot 10^{-9}{[\rm{ls^3}/{m^{2}}]}$, $b_4=125\cdot10^{-6}{[\rm{ls^2}/{m^2}]}$, based on~\cite{Treiber}. 
\begin{table*}[ht!]
	\caption{Performance improvement}
	\label{Table}
	\centering
	\begin{tabular}{|P{1cm}|P{2.2cm}|P{2.2cm}|P{2.2cm}|P{2.2cm}|P{2.2cm}|P{2.4cm}|}
  	\hline
       \multicolumn{1}{|c|}{}& \multicolumn{3}{|c|} {Lyapunov-based/baseline}& \multicolumn{3}{|c|} {RL/baseline} \\
       \hline
       Percent &  backstepping & P & PI & outlet RL & inlet RL & outlet\&inlet RL \\
		\hline
		$J_{\rm TTT}$  & $1.6\%$ & $1.5\%$ & $ 1.5\%$ & $1.4\%$ & $1.6\%$& $1.4\%$    \\
		\hline 
		 $J_{\rm fuel}$ & $3.4\%$ & $2.5\%$ & $ 3.3\%$ & $3.9\%$ & $3.6\%$& $ 4.0\%$   \\
		\hline
     	$J_{\rm comfort}$  &$30.6\%$ & $37.3\%$& $30.1\%$ & $48.3\%$ & $47.8\%$& $47.9\%$     \\
		\hline
	\end{tabular}
\end{table*}
As shown in Table.~\ref{Table}, the performance indices of each controller are compared with their improvement percentage over the baseline setpoint controller. Among the three performance indices, drivers' comfort is the most significantly improved for all Lyapunov-based and RL controls, since the stop-and-go oscillations are suppressed in the closed-loop system. For fuel consumption of total traffic, RL actually performs better than its corresponding Lyapunov-based controller. This may be coincidental, or could be related to RL's control effort and the induced acceleration. In any case, RL provides the flexibility to define the reward as fuel consumption, which can be directly optimized. Backstepping does not provide this capability. We also see that total travel time is only marginally improved, since the traffic dynamics are stabilized to steady state and the average speed remains relatively similar.

\subsection{Comparison study with partial knowledge of system}
In the previous comparative study, we assume a perfect knowledge of the traffic system for the Lyapunov-based controllers. In practice, model parameters are obtained by calibrating the nonlinear ARZ PDE model with field data obtained from loop detectors measuring the traffic flow rate or by high-speed cameras recording vehicle trajectories. The model calibration process can be laborious. More importantly, it is hard to determine some model parameters such as the steady state density, fundamental diagram and relaxation time. These macroscopic models are just that -- model idealizations of reality. In \cite{Yu:auto19}, an adaptive output feedback controller is designed to stabilize the linearized ARZ model with a gradient-based estimator for unknown relaxation time. The adaptive stabilization problem has not been studied when the steady state is uncertain. In traffic field data, it is observed that a certain steady state density in the congested
regime possesses a significant spread of flow rates. Therefore, it is hard to accurately determine the steady state for a real world traffic system that is often periodically evolving, as exhibited in freeway traffic data~\cite{NGSIM}. On the other hand, the steady states in the reward function encode our belief in the current averaged and aggregated traffic condition which could deviate from reality. Moreover, there are measurement errors in raw data. The data needs to be pre-processed leading to more approximation and processing errors. In \cite{Seibold_data} and \cite{Yu:TCST}, the ARZ model was calibrated with the Next Generation Simulation (NGSIM) traffic data which records the trajectories of vehicles on a 500~m freeway segment over a 45~min rush-hour period. Data reconstruction is conducted to obtain the aggregated state values. All in all, the Lyapunov-based control approaches cannot overcome the limitation of model uncertainty when they are implemented in practice. In other words, their performance is limited by the best prediction the PDE model can make. 

Therefore, it is of practical relevance to evaluate the closed-loop performance of the Lyapunov-based and RL controllers given partial or inaccurate models.
Here we investigate the performance of the Lyapunov-based controllers that employ incorrect steady state density, and then compare with RL controllers obtained from a stochastic training process.

We choose two representative scenarios characterized by the steady state density. We define $\rho_r$ as the real steady state density while $\rho^\star =  120 \, \rm{veh}/\rm{km}$ is assumed by the Lyapunov-based feedback controllers. We consider $\rho_r = 115 \, \rm{veh}/\rm{km}$ for the outlet backstepping control case, representing that the actual traffic is lighter than the parameters used by the PDE backstepping controller. In the second scenario, we consider $\rho_r = 125 \, \rm{veh}/\rm{km}$ for the inlet control such that the actual traffic is denser than the value employed by the P controller. For both scenarios, RL controllers are trained in a stochastic environment where a steady state density value of the traffic system $\rho^\star$ is randomly chosen from a uniform distribution of the integer values $\{115, 120, 125\}$ in each episode. Then the RL controllers are validated for each specific scenario.

\begin{figure*}[ht!]
    \centering
    \begin{subfigure}[h]{0.48\textwidth}
        \centering
        \includegraphics[trim = 0.0mm 0.0mm 0.0mm 20.0mm, clip, width=0.90\textwidth]{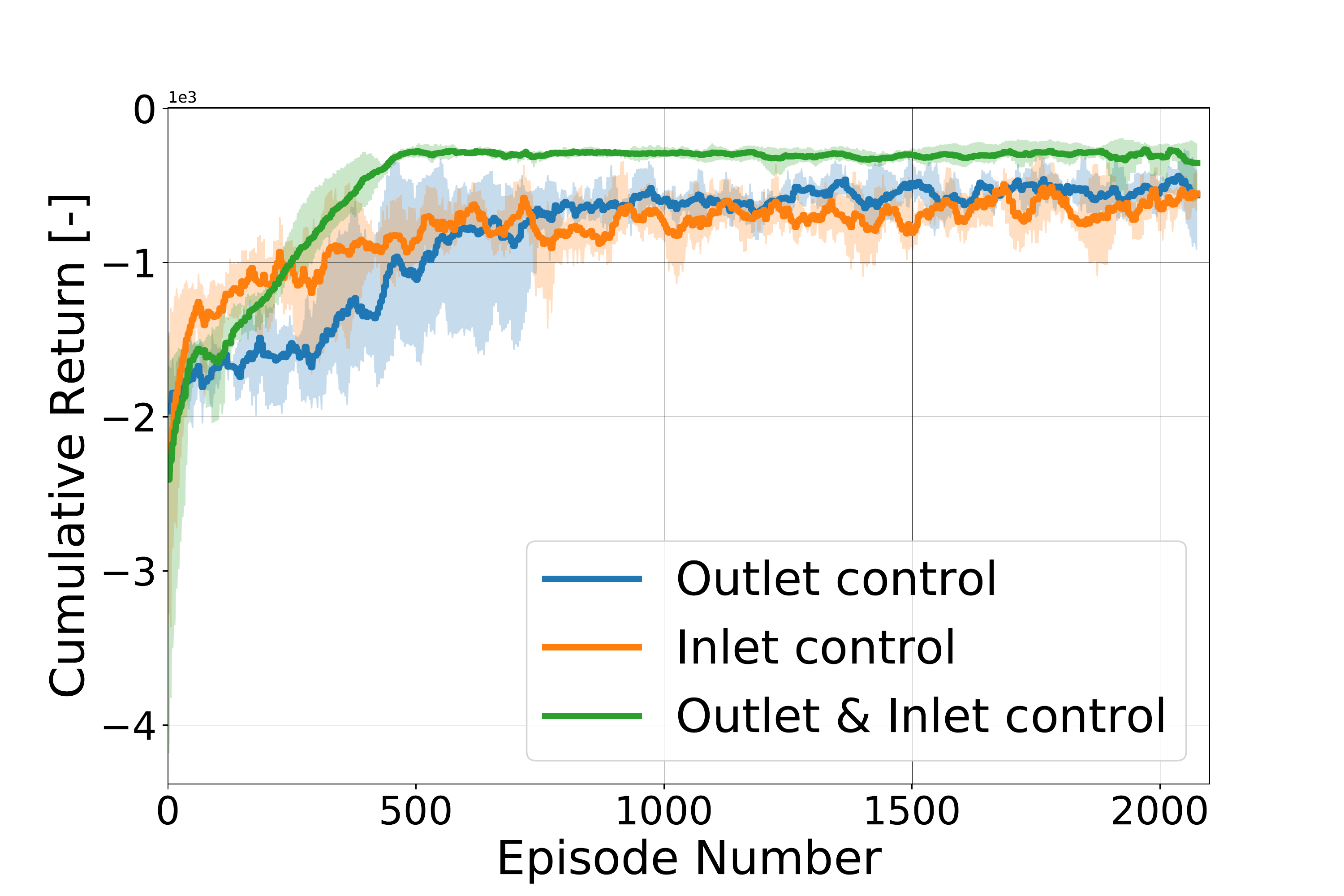}
        \caption{RL learning curve with different initializations of actor-critic network, for randomly generated steady-state values.}
        \vspace{-1ex}
        \label{fig:sto_learning_curve}
    \end{subfigure}
        \hfill
        \begin{subfigure}[h]{0.48\textwidth}
        \centering
        \includegraphics[trim = 0.0mm 0.0mm 0.0mm 0.0mm, clip, width=0.90\textwidth]{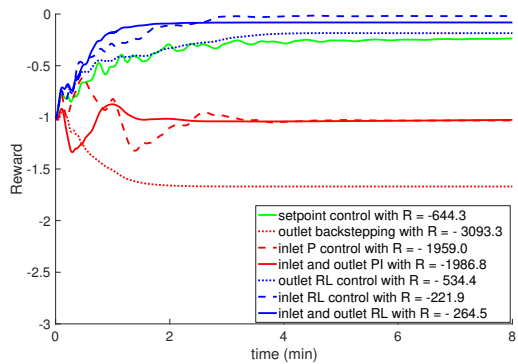}
        \caption{In the system with partial knowledge, the reward evolution of the Lyapunov-based approaches (red), corresponding RL approaches (blue), and setpoint control (green).}
        \vspace{0ex}
        \label{fig:sto_reward}
    \end{subfigure}%
    \caption{The learning curves for RL and reward comparison with the Lyapunov controllers in Scenario 1.}
    \vspace{0ex}
    \label{fig:sto_reward_curves}
\end{figure*}

The RL learning curve is plotted in Fig.~\ref{fig:sto_learning_curve}. Note that the steady state density used to calculate the cumulative reward is randomly chosen at each episode, which affects the initial condition and the associated steady state velocity in \eqref{ref-state} during training. The shaded area represents min/max of different actor-critic parameterizations, while the bold line measures the average performance. Similar to Fig.~\ref{fig:det_reward_curves}, it should be noted that the performance of RL controllers do not monotonically improve in the learning process as the number of training episodes increases. Additionally, cumulative reward convergence for some initiations of the actor-critic parameters is not guaranteed.

We observe that the RL controllers are more adaptable to the stochastic environment.
Figure~\ref{fig:sto_reward} shows all the closed-loop controller results in Scenario 1 where the Lyapunov-based controllers assume a greater steady state density value than the actual traffic environment. It is interesting to find out that all of the RL controllers in blue outperform the setpoint control in green, with larger reward values.  Among them, the inlet RL controller ultimately performs the best. In addition, all the Lyapunov-based feedback controllers in red have smaller rewards, indicating the traffic states do not converge to the actual steady state value after applying the controllers. Assuming an incorrect steady state density value deteriorates the closed-loop stabilization results. The backstepping outlet controller, which performs the best in a system with perfect knowledge from Section.~\ref{det}, turns out to be the worst controller in terms of cumulative reward in Scenario 1, as depicted in Fig.~\ref{fig:sto_reward}. This observation demonstrates that Lyapunov-based control is sensitive to full knowledge of the system dynamics, and even perturbing the steady state values can result in worse performance than a baseline control method.

\subsubsection{Scenario 1 of lighter in-domain traffic} 
\begin{figure*}[ht!]
\hspace*{-0.2cm}
\begin{subfigure}[a]{.6\textwidth}
  \centering
  \includegraphics[width=1.05\linewidth]{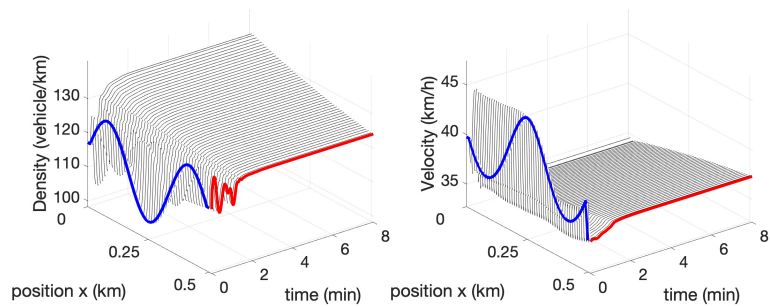}  
 \caption{closed-loop with outlet backstepping controller {\bf in a lighter traffic}}
  \label{fig:bkst_s}
 \hspace*{-0.15cm}
  \includegraphics[width=1.05\linewidth]{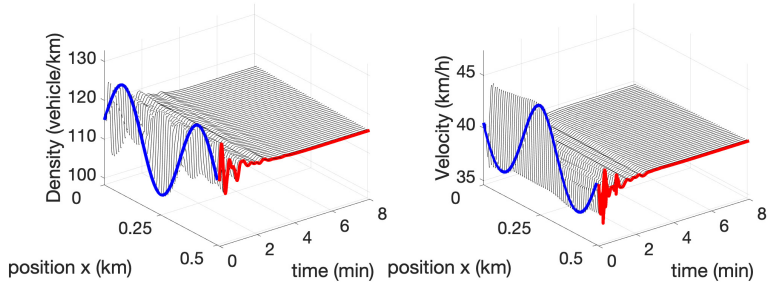}
 \caption{closed-loop with outlet RL controller {\bf in a lighter traffic}}
  \label{fig:RL_outlet_s}
\end{subfigure}
\hspace*{-0.3cm}
\begin{subfigure}[c]{0.4\textwidth}  
\vspace{0.1cm}
       \centering 
            \includegraphics[width=0.9\textwidth]{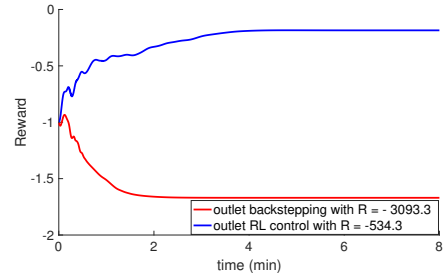}
            \caption{Reward evolution for closed-loop with backstepping and RL controller at outlet \bf in a lighter traffic}
            \label{fig:reward_s_outlet}
\vspace{0.3cm}
            \centering 
            \includegraphics[width=0.9\textwidth]{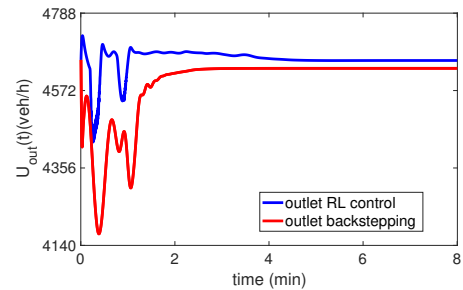}
            \caption{Backstepping and RL control inputs at outlet \bf in a lighter traffic} 
            \label{fig:outlet_s_control}
        \end{subfigure}
\caption{Density and velocity evolution in space and time with control input highlighted with red at outlet }
\label{fig:s_outlet}
\end{figure*}

As shown in Fig. \ref{fig:s_outlet}, we apply the outlet backstepping controller to a nonlinear ARZ model with lighter steady-state traffic density, $\rho_r = 115 \,\rm{veh}/\rm{km}$ than the controller assumes, i.e. $ \rho_r < \rho^\star$. The closed-loop simulation is run for $ 8 \, \rm min$, since the stabilization requires a longer time.  When we apply the PDE backstepping controller in \eqref{bkst_out}, which is constructed using $\rho^\star = 120 \, \rm{veh}/\rm{km}$, the traffic is actually slowed down, leading to a reduction of velocity and an increase of density in Fig. \ref{fig:bkst_s}. The closed-loop traffic is more congested, due to the controller's incorrect model assumptions. In contrast, since the RL controller is trained in an environment with stochastic conditions, where some episodes use a lighter traffic scenario, it successfully stabilizes the traffic state to a uniform value. We find out in Fig.~\ref{fig:RL_outlet_s} that the closed-loop system with the RL outlet controller converges close to the actual steady state. This leads to a larger cumulative reward $R=-534.3$ for RL than the backstepping controller $R=-3093.3$ in Fig.~\ref{fig:reward_s_outlet}. Comparing the backstepping control input and RL outlet control input in Fig. ~\ref{fig:outlet_s_control}, we find out that backstepping applies less control effort than RL. Backstepping's discharging flow rate is less than the actual steady state flow rate in the segment, thus unnecessarily holding vehicles back from leaving the segment and thus exacerbating traffic congestion. As shown in Fig. \ref{fig:bkst_s}, traffic congestion forms and propagates backwards to the segment inlet, whereas there is no such congestion in Fig.~\ref{fig:RL_outlet_s}.

\subsubsection{Scenario 2 of denser in-domain traffic}
\begin{figure*}[ht!]
\hspace*{-0.2cm}
\begin{subfigure}[a]{.6\textwidth}
  \centering
  \includegraphics[width=1.05\linewidth]{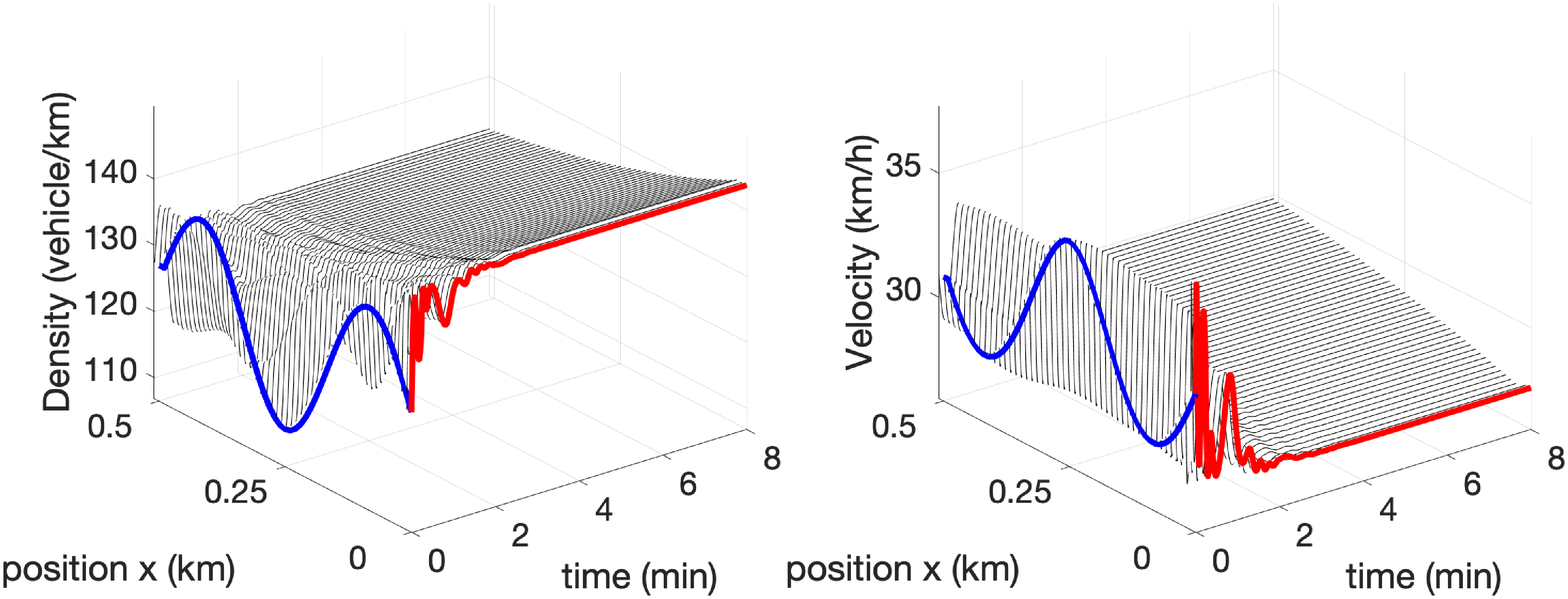}  
 \caption{closed-loop with inlet P controller {\bf in a denser traffic}}
  \label{fig:p_s}
 \hspace*{-0.15cm}
  \includegraphics[width=1.05\linewidth]{RL_p_control.eps}
 \caption{closed-loop with inlet RL controller {\bf in a denser traffic}}
  \label{fig:RL_inlet_s}
\end{subfigure}
\hspace*{-0.3cm}
\begin{subfigure}[c]{0.4\textwidth}  
\vspace{0.1cm}
       \centering 
            \includegraphics[width=0.9\textwidth]{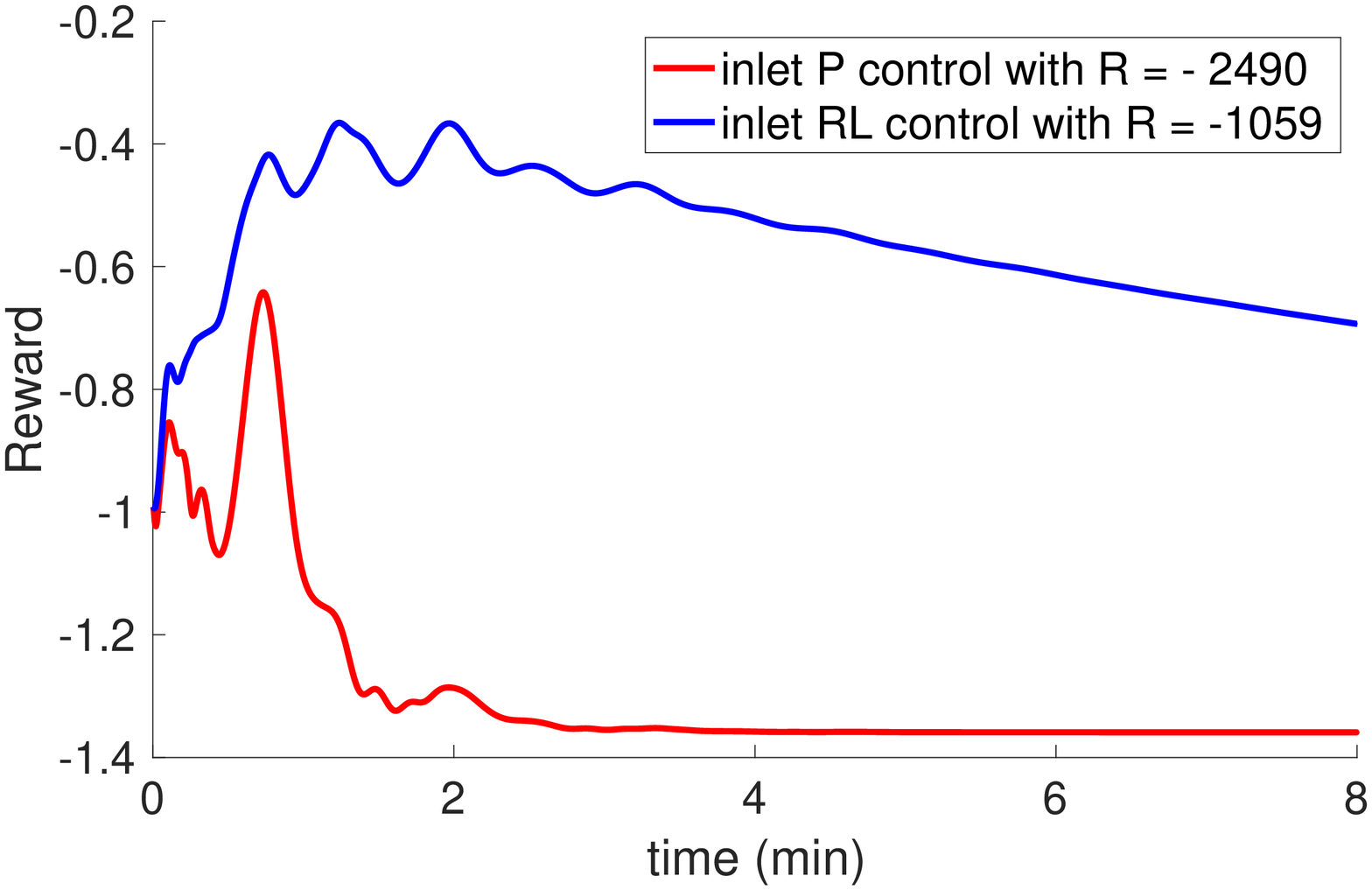}
            \caption{Reward evolution for closed-loop with P and RL controller at inlet \bf in a denser traffic}
            \label{fig:reward_s_inlet}
\vspace{0.3cm}
            \centering 
            \includegraphics[width=0.9\textwidth]{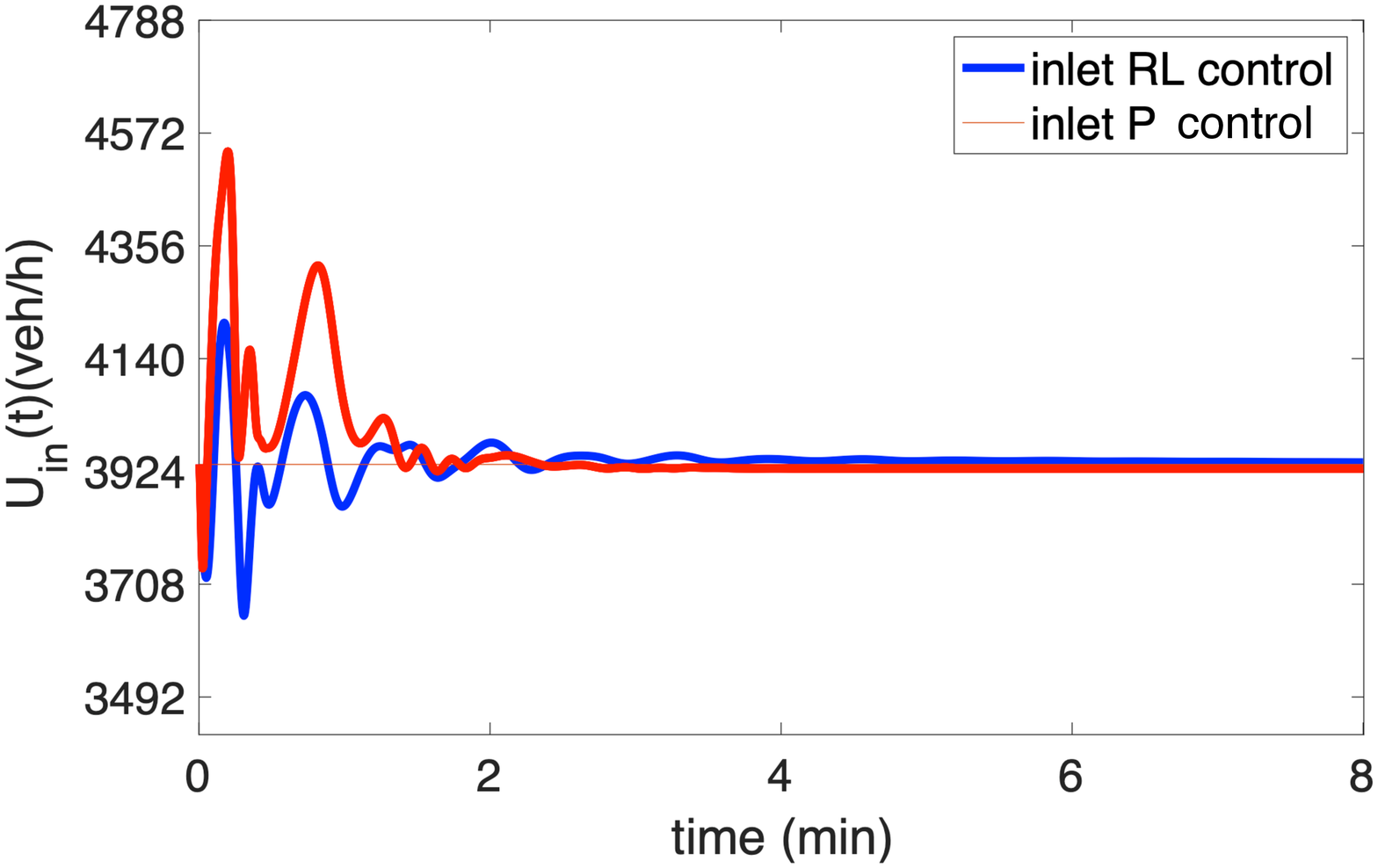}
            \caption{P and RL control inputs at inlet \bf in a denser traffic} 
            \label{fig:inlet_s_control}
        \end{subfigure}
\caption{Density and velocity evolution in space and time with control input highlighted with red at inlet }
\label{fig:s_inlet}
\end{figure*}

In Fig. \ref{fig:s_inlet}, we apply the inlet P controller to a nonlinear ARZ model with denser traffic $\rho_r = 125 \, \rm{veh}/\rm{km}$ than the controller assumes, i.e. $ \rho_r > \rho^\star$. The inlet controller in \eqref{p_in} assumes $\rho^\star = 120\, \rm{veh}/\rm{km}$ to design the control gain. As a result, the P controller applies a larger incoming traffic flow rate than it should, since it receives the feedback of slower velocity than the incorrect reference value. This creates congestion in the downstream traffic. As shown at the inlet of the Fig.~\ref{fig:p_s}, the traffic density increases while velocity reduces. In contrast, the RL controller trained from the stochastic environment alleviates this situation in Fig.~\ref{fig:reward_s_inlet}, and also exhibits an improved cumulative reward. However, the RL controller cannot stabilize the system to the true steady state $\rho_r$, as shown by the blue reward initially increasing and then decreasing. Figure~\ref{fig:inlet_s_control} describes the evolution of the P and RL inlet control inputs. The incoming flow rate is eventually regulated by the inlet controllers to the real steady state value, at the inlet. However, neither the RL controller nor the P controller successfully stabilizes the states to their true steady state values across the segment, as highlighted in red Fig.~\ref{fig:p_s} and Fig.~\ref{fig:reward_s_inlet}. This is also demonstrated by the reward evolution in Fig.~\ref{fig:reward_s_inlet}. However, the RL controller does not diverge as quickly, and in this sense is more robust.

In both scenarios, the performance of the Lyapunov-based controllers deteriorates with small errors to the assumed steady state conditions. Comparing the state and reward evolution with outlet control and perfect system knowledge in Fig.~\ref{fig:bkst},~\ref{fig:reward_bkst}  and Fig.~\ref{fig:bkst_s},~\ref{fig:reward_s_outlet} in Scenario 1, with the same comparison between Fig.~\ref{fig:P},~\ref{fig:reward_p} and Fig.~\ref{fig:p_s},~\ref{fig:reward_s_inlet} in Scenario 2, the Lyapunov-based control systems do not converge to the actual steady states and the corresponding rewards do not converge to 0. On the other hand, the RL controllers trained in a stochastic environment perform worse in the deterministic setting, but are more robust to steady state errors compared to the Lyapunov-based controllers, as shown in Fig.~\ref{fig:s_outlet} for Scenario 1 and in Fig.~\ref{fig:s_inlet}. It should also be noted in  Fig.~\ref{fig:RL_inlet_s} and Fig.~\ref{fig:reward_s_inlet} that RL controllers do not guarantee convergence of the closed-loop system, even though certain improvements are made for reward. 

\section{Conclusion} 


In this work, we address the freeway traffic control problem with a reinforcement learning paradigm. We are motivated by two challenges currently unaddressed in the literature: (i) Controllers natively designed for nonlinear traffic dynamics do not exist, and existing stabilization results are local and derived with linearized systems. (ii) Model-based controllers depend on accurate models. However, even the state-of-art Aw-Rascle-Zhang (ARZ) second-order PDE model does not capture the complex realities of true real-world traffic flow, making model calibration a significant challenge. Even if one does successfully calibrate a ARZ model on a dataset, the parameters can shift over time due to weather, evolving commute patterns, discrete events, and more. RL methods provide an opportunity to address this practical challenge.

Ramp metering boundary controllers are developed to stabilize the freeway traffic oscillations. The second-order ARZ PDE model is employed to describe the traffic dynamics. We formulate the control of the hyperbolic PDE model as a reinforcement learning problem, where the control objective is encoded into the reward function. RL controllers are designed for various boundary control schemes, along with corresponding Lyapunov-based feedback control strategies. The nonlinear ARZ PDE model is used as the environment to simulate and generate data for reinforcement learning. Among the RL algorithms, the actor-critic scheme is considered, and specifically the PPO algorithm has been adopted. The closed-loop performance of RL controllers is compared against Lyapunov-based state-feedback controllers, including PDE backstepping, P, and PI controllers. 

Simulation results show that RL controllers achieve comparable performance with the Lyapunov-based controllers in a traffic system with perfect knowledge of model parameters, with slightly longer convergence times. The analysis becomes more interesting when we consider the reality of model uncertainty and imperfect knowledge.
Specifically, the RL controllers outperform the Lyapunov-based controllers when the assumed steady-state conditions contain errors. Importantly, we find the Lyapunov-based controller loses its stabilizing properties, whereas the RL controllers obtained from a stochastic training process maintain stability or at least slow down divergence. Nevertheless, the RL controllers are obtained from thousands of iterative training episodes, and reward convergence is not guaranteed. For practical considerations, RL controllers should be developed in a simulation environment where trial and error are allowed. 


{ One of the main advantages of the proposed RL formulation is that it does not require any prior knowledge of the traffic model, besides a simulator, but can achieve comparative and possibly better closed-loop performance compared to a Lyapunov-based approach. Future work may validate the RL-based controller on a freeway network comprised of several segments, merging and diverging lanes. Another direction is to utilize the Lyapunov-based approach with stability guarantees to improve the RL controllers to obtain faster learning rates and more robust performance.}

\section*{Acknowledgement}
We would like to thank Eugene Vinitsky and Aboudy Kreidieh from the Mobile Sensing Lab at the University of California, Berkeley for their insights and help at the beginning of this work.


\begin{thebibliography}{99}  


 

\bibitem{IArel}
I. Arel, C. Liu, T. Urbanik, and A.G. Kohls, 
\newblock {``Reinforcement learning-based multi-agent system for network traffic signal control,"} \newblock {\em IET Intelligent Transport Systems}, 4(2), pp.128-135, 2010.

\bibitem{AW} 
A. Aw,  and M. Rascle,  
\newblock {``Resurrection of" second order" models of traffic flow."} 
\newblock  {\em SIAM journal on applied mathematics,} 60(3), 916-938, 2000.

	
\bibitem{Bekiaris-Liberis:18}
N. Bekiaris-Liberis, and M. Krstic, 
\newblock ``Compensation of actuator dynamics governed by quasilinear hyperbolic PDEs,"
\newblock  {\em Automatica}, vol.92, pp.29-40, 2018.

\bibitem{bell} 
 F. Belletti,  M. Huo, X. Litrico and  A. M. Bayen,  
 \newblock {``Prediction of traffic convective instability with spectral analysis of the Aw–Rascle–Zhang model."}  
 {\em Physics Letters A}, 379(38), 2319-2330, 2015.



\bibitem{BELLETTI} 
F. Belletti,  D. Haziza,  G.Gomes, and A.M. Bayen,  
\newblock {``Expert level control of ramp metering based on multi-task deep reinforcement learning."} 
\newblock {\em IEEE Transactions on Intelligent Transportation Systems }, pp. 1198-1207, 2017.


 
\bibitem{bertsekas2005dynamic}
D. Bertsekas, 
\newblock {``Dynamic programming and optimal control."} 
\newblock {\em Athena scientific Belmont}, 2005.


 
 \bibitem{Daganzo:11}
C. F. Daganzo, 
 \newblock {``On the macroscopic stability of freeway traffic."}   
 {\em Transportation Research Part B: Methodological}, 45(5), 782-788, 2011.
 
 
 
\bibitem{MARLIN} 
S. El-Tantawy, B. Abdulhai, and H. Abdelgawad,
\newblock {``Multiagent reinforcement learning for integrated network of adaptive traffic signal controllers (MARLIN-ATSC): methodology and large-scale application on downtown Toronto,"}  \newblock {\em IEEE Transactions on Intelligent Transportation Systems} vol.14.3, pp. 1140-1150, 2013. 
 
 \bibitem{Seibold_data} 
 S. Fan, and  B.  Seibold,
\newblock {``Data-fitted first-order traffic models and their second-order generalizations: Comparison by trajectory and sensor data."}
\newblock  {\em  Transportation research record,} 2391(1), 32-43, 2013. 
 
\bibitem{Seibold1}
 M. R. Flynn, A. R. Kasimov,  J. C.Nave,  R. R. Rosales, B.  Seibold, 
 \newblock {``Self-sustained nonlinear waves in traffic flow."} 
  {\em Physical Review E}, 79(5), 056113, 2009.
  
\bibitem{Farahmand17} 
A. M. Farahmand,  S. Nabi,and  D.N. Nikovski, 
\newblock {``Deep reinforcement learning for partial differential equation control,"}\newblock {\em In 2017 American Control Conference (ACC)}, pp. 3120-3127, 2017, May.

\bibitem{Farahmand18}
 Y. Pan, A.M. Farahmand, M. White, S.Nabi, P. Grover, and D. Nikovski, 
 \newblock {`` Reinforcement learning with function-valued action spaces for partial differential equation control,"} arXiv preprint arXiv:1806.06931, 2018.


\bibitem{NGSIM}
FHWA, U.S. Department of Transportation. Next Generation Simulation
(NGSIM). http://ops.fhwa.dot.gov/trafficanalysistools/ngsim.html.


\bibitem{Frihauf}
  P. Frihauf, M. Krstic, and T. Basar, ``Finite-horizon LQ control for unknown discrete-time linear systems via extremum seeking,” {\em  European Journal of Control}, vol. 19, pp. 399-407, 2013.


\bibitem{Gomes:06}
 G. Gomes and  R. Horowitz, 
 \newblock {``Optimal freeway ramp metering using the asymmetric cell transmission model."} 
\newblock  {\em   Transportation Research Part C: Emerging Technologies,}  14(4),244-262, 2016. 

\bibitem{JHan18}
J. Han,  A. Jentzen, and  E. Weinan, \newblock {``Solving high-dimensional partial differential equations using deep learning."}  \newblock {\em Proceedings of the National Academy of Sciences}, 115(34), pp.8505-8510, 2018.


\bibitem{Kakade} 
S. Kakade and J. Langford, 
\newblock {``Approximately optimal approximate reinforcement learning."} \newblock {\em In 2002 International Conference on Machine Learning (ICML)}, pp. 267-274. IEEE, 2002.

\bibitem{Karafyllis:19} 
 I. Karafyllis, and  M. Papageorgiou, 
\newblock {``Feedback control of scalar conservation laws with application to density control in freeways by means of variable speed limits."} 
\newblock {\em Automatica}, no.105, pp.228-236, 2019.

\bibitem{Karafyllis:18} 
I. Karafyllis,  N. Bekiaris-Liberis, and M. Papageorgiou, 
\newblock {``Feedback Control of Nonlinear Hyperbolic PDE Systems Inspired by Traffic Flow Models,"}
\newblock  {\em IEEE Transactions on Automatic Control,} 2018.


\bibitem{intro_bkst}
 M. Krstic, and  A. Smyshlyaev.
 \newblock {``Boundary control of PDEs: A course on backstepping designs."} \newblock {\em Society for Industrial and Applied Mathematics}, 2008.


\bibitem{Nesic}
 S. Z. Khong, D. Nesic, and M. Krstic, 
 ``A non-model based extremum seeking approach to iterative learning control,” 
 {\em  Automatica}, vol. 66, pp. 238-245, 2016.

\bibitem{killingsworth}
 N. Killingsworth and M. Krstic,
 ``PID tuning using extremum seeking,” 
{\em  Control Systems Magazine}, vol. 26, pp. 70-79, February 2006.

 
\bibitem{RL_Pytorch} 
I. Kostrikov, 
\newblock {``PyTorch Implementations of Reinforcement Learning Algorithms,"}
In Github Repository, https://github.com/p-christ/Deep-Reinforcement-Learning-Algorithms-with-PyTorch.



\bibitem{LaxWendroff} 
P. Lax,  and  B. Wendroff. 
\newblock {``Systems of conservation laws."} 
\newblock {\em Communications on Pure and Applied mathematics}, 13(2), pp. 217--237, 1960.


\bibitem{LiLi}
L. Li,  Y. Lv, and  F.Y.Wang, 
\newblock {``Traffic signal timing via deep reinforcement learning,"} \newblock {\em IEEE/CAA Journal of Automatica Sinica}, 3(3), pp.247-254, 2016.

 
\bibitem{RLlib} 
E. Liang, R. Liaw, P. Moritz,  R. Nishihara, R. Fox, K.Goldberg, J. Gonzalez, M. Jordan, and I. Stoica,   
\newblock {``RLlib: Abstractions for distributed reinforcement learning."} \newblock {\em In 2018 International Conference on Machine Learning (ICML)}, preprint arXiv:1712.09381, IEEE, 2018.

\bibitem{LW} 
M. J. Lighthill, and  G. B. Whitham, 
\newblock {``On kinematic waves: II. A theory of traffic flow on long crowded roads."}  
\newblock{\em Proceedings of the Royal Society of London. Series A, Mathematical and
Physical Sciences (1934–1990),} vol.229 (1178), pp. 317–345, 1955.



\bibitem{Papageorgiou:90} 
M. Papageorgiou, J. M. Blosseville, H. Haj-Salem,
\newblock {``Modelling and real-time control of traffic flow on the southern part of Boulevard P{\'e}riph{\'e}rique in Paris – Part I: modelling,"}  
\newblock {\em Transportation Research Part A}, vol.24 (4), pp.345–359, 1990. 

\bibitem{Radenkovic}
 M. Radenkovic and M. Krstic,  ``Adaptive control via extremum seeking: Global stabilization and consistency of parameter estimates,”  {\em  IEEE Transactions on Automatic Control}, vol. 62, pp. 2350-2359, 2017.
 
\bibitem{R}
P. I. Richards, 
\newblock {``Shock waves on the highway."} 
\newblock  {\em Operations Research,}  vol.4 (1), pp.42–51, 1956.



 \bibitem{TRPO} 
J. Schulman, S. Levine, P. Abbeel, M.Jordan and P. Moritz, 
\newblock {``Trust Region Policy Optimization."} 
\newblock {\em In 2015 International Conference on Machine Learning (ICML)}, pp. 1889-1897. IEEE, 2015.

\bibitem{PPO} 
J. Schulman, F. Wolski, P. Dhariwal, A. Radford, O.Klimov, 
\newblock {``Proximal policy optimization algorithms."} arXiv preprint arXiv:1707.06347, 2017.


\bibitem{Seibold2}
B. Seibold, M. R. Flynn,  A. R. Kasimov, R. R. Rosales,  
\newblock {``Constructing set-valued fundamental diagrams from jamiton solutions in second order traffic models."}  arXiv preprint arXiv:1204.5510, 2012.


\bibitem{sutton2000policy}
R. Sutton, D. McAllester, s. Singh, and Y.Mansour, 
\newblock {``Policy gradient methods for reinforcement learning with function approximation." }
\newblock {\em In Advances in neural information processing systems}. pp. 1057 -1063, 2000.



\bibitem{Treiber}
 M. Treiber, and  A. Kesting, 
\newblock {``Traffic Flow Dynamics: Data, Models and Simulation,"}  
\newblock  {\em Springer-Verlag Berlin Heidelberg}, 2013.

\bibitem{Vinitsky18}
E. Vinitsky, K. Parvate, A. Kreidieh, C. Wu, and  A. Bayen, 
\newblock {``Lagrangian control through deep-rl: Applications to bottleneck decongestion."} \newblock {\em In 2018 21st International Conference on Intelligent Transportation Systems (ITSC)}, pp. 759-765. IEEE, 2018, November.

\bibitem{YWang19}
Y. Wang, Z. Shen, Z. Long, and B. Dong, \newblock {``Learning to discretize: solving 1D scalar conservation laws via deep reinforcement learning"}. arXiv preprint arXiv:1905.11079 , 2019.



\bibitem{Wu17}
C. Wu, A. Kreidieh, E. Vinitsky, and A.M. Bayen, 
\newblock {``Emergent behaviors in mixed-autonomy traffic."}  In Conference on Robot Learning, pp. 398-407. 2017.



\bibitem{Yu:auto19} 
H. Yu, and M. Krstic,
\newblock ``Traffic congestion control for Aw–Rascle–Zhang model,” \newblock  {\em Automatica}, vol. 100, pp. 38-51, 2019

\bibitem{Huan2}  
H. Yu, and  M. Krstic,  \newblock {``Varying Speed Limit Control of Aw-Rascle-Zhang Traffic Model."} \newblock {\em In 2018 21st International Conference on Intelligent Transportation Systems (ITSC)} , pp. 1846-1851. IEEE, 2018.


\bibitem{Yu:TCST}
H. Yu,  Q. Gan,  A. Bayen, and  M. Krstic,
\newblock ``PDE Traffic Observer Validated on Freeway Data," \newblock  {\em  IEEE Transactions on Control Systems Technology}. http://dx.doi.org/10.1109/TCST.2020.2989101, 2020.


\bibitem{Yu:TAC}
H. Yu, M. Diagne,  L. Zhang, and  M. Krstic, 
\newblock ``Bilateral boundary control of moving shockwave in LWR model of congested traffic, "
\newblock  {\em IEEE Transactions on Automatic Control,} 2020.

\bibitem{Yu:JDSMC}
H. Yu, S. Koga, T. R. Oliveira, and M. Krstic, 
\newblock ``Extremum seeking for traffic congestion control with a downstream bottleneck," 
\newblock  {\em ASME Journal of Dynamic Systems, Measurement, and Control},  DOI: 10.1115/1.4048781, 2020.






\bibitem{Zhang} 
H. M. Zhang,
\newblock {``A non-equilibrium traffic model devoid of gas-like behavior."}  
\newblock  {\em  Transportation Research Part B: Methodological,} 36(3), 275-290, 2002.

\bibitem{Liguo:19}
L. Zhang, C. Prieur, and J. Qiao.  
\newblock {``PI boundary control of linear hyperbolic balance laws with stabilization of ARZ traffic flow models."} \newblock {\em  Systems \& Control Letters}, 123, 85-91, 2019.


\end{thebibliography}
\end{document}